\documentclass[11pt]{article}
\usepackage{amssymb, amsmath, amscd, amsthm}
\usepackage{mathrsfs}
\usepackage[all]{xy}

\usepackage{amssymb, amsmath, amscd, amsthm}
\usepackage[OT2,OT1]{fontenc}
\usepackage[all]{xy}

       {\arraycolsep 0.14em\begin{eqnarray}}{\end{eqnarray}}
\newenvironment{Eqnarray*}%
       {\arraycolsep 0.14em\begin{eqnarray*}}{\end{eqnarray*}}
       {\arraycolsep 0.14em\begin{array}}{\end{array}}

\theoremstyle{plain}
\newtheorem{Thm}{Theorem}[section]
\newtheorem{Lem}[Thm]{Lemma}
\newtheorem{Prop}[Thm]{Proposition}
\newtheorem{Cor}[Thm]{Corollary}

\theoremstyle{definition}
\newtheorem{Def}[Thm]{Definition}

\theoremstyle{remark}
\newtheorem{Rem}[Thm]{Remark}

\numberwithin{equation}{section}

\def\action{\, {\scriptscriptstyle \stackrel{\circ}{{}}} \, }

\def\C{\mathbb C}  \def\bD{\mathbb D}  \def\bH{\mathbb H} 
  \def\bP{\mathbb P}  
\def\R{\mathbb R}  \def\bS{\mathbb S}  \def\Z{\mathbb Z}  
\def\bV{\mathbb V}

\def\a{\alpha}  \def\b{\beta}     \def\g{\gamma}   \def\d{\delta}
\def\c{\theta}  \def\l{\lambda}   \def\p{\phi}     \def\s{\sigma}
\def\t{\tau}    \def\vp{\varphi}       \def\y{\eta}
\def\z{\zeta}   \def\w{\omega}    \def\e{\varepsilon}
\def\x{\xi}       \def\k{\kappa}

   \def\cD{\mathcal D}   \def\cF{\mathcal F}
      
   \def\cM{\mathcal M}  \def\cN{\mathcal N}
     
\def\cQ{\mathcal Q}  \def\cR{\mathcal R}   
      
\def\cW{\mathcal W}   \def\cZ{\mathcal Z}  

\def\scrH{\mathscr H}

\DeclareMathAlphabet{\mathzap}{OT1}{pzc}{m}{it}

\def\cyr{%
\renewcommand\rmdefault{wncyr}%
\renewcommand\sfdefault{wncyss}%
\renewcommand\encodingdefault{OT2}%
\normalfont
\selectfont}
\DeclareTextFontCommand{\textcyr}{\cyr}

\newcommand{\pye}{\mbox{\cyr p}}

\newcommand{\xye}{\mbox{\cyr x}}

\def\Ga{\mathfrak a}
\def\Gf{\mathfrak f}       
  \def\Gm{\mathfrak m}

\def\GC{\mathfrak C}  \def\GD{\mathfrak D}  
\def\GS{\mathfrak S}

\def\U{\operatorname{U}}

\def\SO{\operatorname{SO}}        \def\Sp{\operatorname{Sp}}

\def\Aut{\operatorname{Aut}}

\def\Span{\operatorname{Span}}

\def\Imag{\operatorname{Im}}    \def\Real{\operatorname{Re}}      

\def\RP{\R\bP}
\def\CP{\C\bP} 

\def\del{\partial}

\def\({ \left( }     \def\){ \right) }
\def\<{ \left\langle } \def\>{ \right\rangle }
\def\hasira{\rule{0mm}{2eX}}

\def\ul{\underline}

\newcommand{\pd}[1]{\frac{\del}{\del #1}}
\newcommand{\PD}[2]{\frac{\del #1}{\del #2}}
\newcommand{\inn}[2]{\< #1 ,#2 \>}

\title{Wave equations and the LeBrun-Mason correspondence}
\author{Fuminori Nakata\thanks{This work 
 is partially supported by Research Fellowships of the Japan Society for the 
 Promotion of Science of Young Scientists. }}

\date{}

\begin{document}
\maketitle

\abstract{ 
The LeBrun-Mason twistor correspondences for $S^1$-invariant 
self-dual Zollfrei metrics are explicitly established. 
We give explicit formulas for the general solutions of 
the wave equation and the monopole equation on the 
de Sitter three-space under the assumption for the tameness at infinity 
by using Radon-type integral transforms, 
and the above twistor correspondence is 
described by using these formulas. 
We also obtain a critical condition for the LeBrun-Mason twistor spaces, 
and show that the twistor theory does not work well for 
twistor spaces which do not satisfy this condition. }

\vspace{8mm}
{\noindent
{{\it Mathematics Subject Classifications} (2000) : 
  53C28,  
  35L05,  
  53C50,  
  32G10. \\  
  {\it Keywords} : \  twistor method, holomorphic disks, 
   indefinite metric, wave equation, \\ 
   \hspace{17mm} monopole equation, Radon transform.}

\section{Introduction}
The twistor theory concerning holomorphic disks, developed by 
C.~LeBrun and L.~J.~Mason, is now progressing steadily 
(see \cite{bib:LeBrun06,bib:LM02,bib:LM05,bib:LM08,
bib:Nakata07I,bib:Nakata07II,bib:Nakata09}). 
In general, LeBrun-Mason type twistor correspondence 
is characterized in the following way: 
\begin{itemize}
 \item the twistor space is given by the pair of a complex manifold $Z$ and 
   a totally real submanifold $P$ in $Z$, 
 \item corresponding objects to the `twistor lines' in the ordinary twistor theory 
   \cite{bib:AHS,bib:Hitchin82,bib:Penrose76} 
   are given by the {\it holomorphic disks} on $Z$ with boundaries lying on $P$, 
 \item a natural differential geometric structure is induced on 
    the parameter space $M$ of the family of holomorphic disks, 
 \item the induced structure is of low regularity in general, and satisfies some 
    {\it global conditions} which give a strong restriction on the topology on $M$, and 
 \item conversely, the twistor space $(Z,P)$ is obtained from 
   such differential geometric structure. 
\end{itemize}
In this article, we mainly deal with the non-rigid case of 
the LeBrun-Mason correspondence for self-dual 
conformal structures \cite{bib:LM05}. 
In this case, the twistor space is a pair $(\CP^3,P)$ where $P$ is an 
embedded $\RP^3$ sufficiently close to the standard one, 
and the corresponding geometry is a self-dual indefinite conformal structure $[g]$ 
on $S^2\times S^2$ of signature $(--++)$. 
In this case, the required global condition for $[g]$ is the {\it Zollfrei condition}, 
that is, every maximal null geodesic of $[g]$ is closed (cf.\cite{bib:Guillemin89}). 
In \cite{bib:LM05}, it is shown that any self-dual indefinite conformal structure on 
$S^2\times S^2$ sufficiently close to the standard one is automatically Zollfrei, 
and that such conformal structures one-to-one corresponds with the twistor 
spaces $(\CP^3,P)$ in the above sense. 

On the other hand, before LeBrun and Mason develop the above theory, 
infinitely many examples of self-dual indefinite metrics on 
$S^2\times S^2$ are obtained by K.~P.~Tod \cite{bib:Tod93}, 
and independently by H.~Kamada \cite{bib:Kamada05}. 
Tod constructed $S^1$-invariant self-dual indefinite 
metrics on $S^2\times S^2$ via method analogous to 
what is called LeBrun's hyperbolic ansatz \cite{bib:LeBrun91}. 
Kamada investigated compact scalar-flat indefinite K\"ahler surfaces with 
Hamiltonian $S^1$-symmetry. It is known that such surface is automatically 
self-dual, and Kamada proved that 
such structure is admitted only on $\CP^1\times\CP^1$. 
Kamada also constructed infinitely many examples of such structures 
containing Tod's examples. 
Since Tod's and Kamada's examples contain the self-dual metrics sufficiently 
close to the standard one, 
at least some of them must be Zollfrei 
by the above results by LeBrun and Mason. 
So the natural question is the following: 
\begin{itemize}
 \item Are the metrics constructed by Tod or Kamada all Zollfrei? 
 \item If they are Zollfrei, can we establish the LeBrun-Mason correspondences 
   for them? 
\end{itemize}
We show that these problems are settled positively, 
which is the main theorem in this article (Theorem \ref{Thm:correspondence}).

To attack the above problems, 
we first study the {\it wave equation} on the 
three-dimensional {\it de Sitter space} $S^3_1$ in 
Section \ref{Sect:Wave} and \ref{Sect:Inverse}. 
We introduce Radon-type integral transforms, 
and show that any solution of the wave equation on $S^3_1$ 
which is tame at infinity is obtained from a function on $S^2$ by 
applying these transforms (Theorem \ref{Thm:Inverse}). 
As a consequence, we see that any solution of the wave equation on $S^3_1$ 
which is tame at infinity carries a symmetry 
which we call the {\it oddness}.  

We next study the monopole equation on $S^3_1$ 
in Section \ref{Sect:Monopole_equation}. 
We introduce the notion of a {\it monopole potential} 
and show that any gauge equivalent class of monopole solutions 
one-to-one corresponds with a monopole potential. 
Further, based on the above results for the wave equation, 
we show that gauge equivalent classes of monopole solutions 
which are tame at infinity one-to-one correspond with  
functions on $S^2$ which we call {\it generating functions} 
(Theorem \ref{Thm:3objects}).

If we follow Kamada's formulation, we can construct self-dual indefinite metrics 
on $S^2\times S^2$ from monopole solutions on $S^3_1$ 
satisfying some extra conditions. 
In light of this construction, we introduce the notion of 
{\it admissible monopoles} by which we obtain the 
self-dual metrics on $S^2\times S^2$ via Kamada's construction. 
By the results above, we see that the admissible monopoles 
are obtained from generating functions satisfying certain condition 
corresponding to the admissibility. 
We remark that by this method we obtain all the monopole solutions 
by which Kamada's construction works. In particular, our method covers 
all the examples obtained by Tod and Kamada.

In the latter half of this article 
(Section \ref{Sect:Local_theory}, \ref{Sect:Standard_model} and 
  \ref{Sect:Twistor_correspondence}), 
we establish the LeBrun-Mason correspondence 
for the above obtained self-dual metrics on $S^2\times S^2$. 
We set a twistor space $(\CP^3,P_h)$ 
for each generating function $h$ on $S^2$, 
and we show that we can establish the LeBrun-Mason correspondence between 
the twistor space $(\CP^3,P_h)$ and the self-dual metric on $S^2\times S^2$ 
obtained from the monopole solution corresponding to $h$ 
if the monopole is admissible. 
In particular we see that the self-dual metrics on $S^2\times S^2$ 
obtained from admissible monopoles are all Zollfrei. 
We also study the non-admissible case, and 
show that the twistor space $(\CP^3,P_h)$ carries an unexpected 
property for holomorphic disks in this case 
(Proposition \ref{Prop:non-admissible_case}).

The results in this article is also considered as the LeBrun-Mason theory version 
of the Jones-Tod reduction theory \cite{bib:JT}. 
In contrast, in \cite{bib:Nakata07I,bib:Nakata07II}, the author studied 
the LeBrun-Mason theory version of the 
Dunajski-West reduction theory \cite{bib:Calderbank,bib:DW}. 
Particularly in \cite{bib:Nakata07I}, 
we obtain infinitely many self-dual indefinite Zollfrei conformal structures on 
$S^2\times S^2$ {\it with singularity}, 
and their LeBrun-Mason correspondences are established 
by making use of the {\it Radon transform} on $\R^2$. 
Though it seems that there are no direct relation between 
these previous works and the results in this article, 
these results seem to insist the significance of the Radon transform 
as a tool in the study of LeBrun-Mason theory.

\section{Wave equation on the de Sitter 3-space}
\label{Sect:Wave}

In this section, we introduce the wave equation on the de Sitter 3-space. 
Then we introduce integral transforms and show that 
we can get solutions of the wave equation by these transforms.

\paragraph{The space of small circles}
Let $\bS^2=\{(u_1,u_2,u_3)\in \R^3 \mid u_1^2+u_2^2+u_3^2=1\}$ 
be the unit sphere equipped with the standard metric and 
$(S^3_1,g_{S^3_1})$ be the de Sitter 3-space defined by 
$$ \begin{aligned} 
 S^3_1 &:=\{ (x_0,x_1,x_2,x_3)\in \R^4 \mid 
 -x_0^2+x_1^2+x_2^2+x_3^2=1 \}, \\
 g_{S^3_1} &:= (-dx_0^2+dx_1^2+dx_2^2+dx_3^2)|_{S^3_1}. 
 \end{aligned} $$ 
We identify $S^3_1$ with $\R\times \bS^2$ via the diffeomorphism 
$\R\times\bS^2 \overset\sim \longrightarrow S^3_1$ given by 
$$ (t,y) \longmapsto (x_0,x_1,x_2,x_3) = 
 (\sinh t,\cosh t \, (y_1,y_2,y_3)). $$

For each $(t,y)\in \R\times\bS^2 \cong S^3_1$, we define 
$$ \Omega_{(t,y)} := \{ u\in \bS^2 \mid u\cdot y > \tanh t \} $$ 
which is an open set on $\bS^2$ bounded by a small circle. 
By the correspondence $(t,y) \leftrightarrow \del\Omega_{(t,y)}$, 
we identify the de Sitter space $S^3_1$ with the space of oriented 
small circles in $\bS^2$. 
Notice that the subset 
$S^2_o:=\{(t,y)\in S^3_1\mid  t=0\}$, which is called the 
{\em neck sphere}, 
corresponds to the space of big circles on $\bS^2$. 
Now let us fix $(t,y)\in S^3_1$ 
and take vectors $y^\perp_1, y^\perp_2 \in \bS^2$  
so that $\{y^\perp_1, y^\perp_2, y\}$ 
gives a basis of $\R^3$ with compatible orientation. 
We define a map $\g_{(t,y)} : S^1 \to \bS^2$ by 
\begin{equation} \label{eq:parameter_of_circle}
 \g_{(t,y)}(\p)= \frac{\cos\p}{\cosh t}\, y^\perp_1 
 + \frac{\sin\p}{\cosh t}\, y^\perp_2 + \tanh t \ y , 
\end{equation}
which gives an oriented parametrization of 
the small circle $\del\Omega_{(t,y)}$. 
We will see later in Section \ref{Sect:Standard_model} that 
the above identification between $S^3_1$ and the space of small circles on 
$\bS^2$ is naturally arisen from the LeBrun-Mason correspondence for 
Einstein-Weyl 3-fold.

\paragraph{The wave equation}
The {\it wave equation} on the 
de Sitter space $S^3_1$ is given by 
\begin{equation} \label{eq:wave_eq_def}
 \Box V:= *d*d V= 0 
\end{equation} 
where $V$ is a smooth function on $S^3_1$ and
 $*$ is the Hodge's operator on $S^3_1$ with respect to 
the indefinite metric $g_{S^3_1}$ and the natural orientation on 
$S^3_1\simeq \R\times\bS^2$

We fix the notations for the operators on $\bS^2$ as follows: 
let $\check *$ be the Hodge's operator, 
$\check d$ be the exterior derivative and 
$\Delta_{\bS^2}$ be the Laplace operator. 
We also use the same notations $\check *$, $\check d$ and $\Delta_{\bS^2}$ 
for the fiberwise operators on $S^3_1$ as a 
$\bS^2$-bundle $S^3_1=\R\times\bS^2\to \R$. 
For any 1-form $\y$ on $S^3_1$ satisfying $\y(\del_t)=0$ 
(where $\del_t=\pd{t}$), 
we obtain 
\begin{equation} 
 *\y = - dt \wedge (\check *\y) \qquad \text{and} \qquad  
 *(dt\wedge\y) = - * \y. 
\end{equation}
If we denote the volume form on $\bS^2$ by $\w_{\bS^2}$, then 
we have $*\, dt=-\cosh^2\! t\,\, \w_{\bS^2}$. 
For a smooth function $V$ on $S^3_1$, 
\begin{equation} \label{eq:*dV}
 \begin{aligned} 
 d \, V &= V_tdt + \check d\, V, \\ 
 *\, d\, V &= 
   - V_t\cosh^2\! t\,\, \w_{\bS^2} - dt \wedge (\check *\,\check d\, V), \\ 
 d*d\, V &= -(V_t\cosh^2\! t)_t\, dt\wedge \w_{\bS^2} 
    + dt \wedge \check d \, \check *\, \check d\, V, 
 \end{aligned}
\end{equation}
where $V_t=\del_tV$ and so on. 
Hence the wave equation \eqref{eq:wave_eq_def} is written as 
\begin{equation} \label{eq:WE}
  \( -\frac{\del^2}{\del t^2} - 2\tanh t \, \frac{\del}{\del t} + 
   (\cosh t)^{-2}\Delta_{\bS^2} \) V =0. 
\end{equation}

\paragraph{Function spaces}
Let us denote the antipodal map on $\bS^2$ by $\a$. 
We also define an involution on $S^3_1$ by 
$ \sigma : (t,y) \mapsto (-t,-y). $
If we identify $S^3_1$ with the space of oriented small circles on $\bS^2$, 
$\s$ corresponds to the orientation reversing operation 
for each oriented small circle. 
Let us denote by $C^\infty(\bS^2)$ and by $C^\infty(S^3_1)$ the 
space of real valued smooth functions on $\bS^2$ and on $S^3_1$ respectively. 
We set 
 $$ \begin{aligned} 
   C^\infty_{\text{even}}(\bS^2) &:= 
    \{ h\in C^\infty(\bS^2) \mid h=h\action \a \}, \\ 
   C^\infty_{\text{odd}}(\bS^2) &:= 
    \{ h\in C^\infty(\bS^2) \mid h=-h\action \a \}, 
  \end{aligned} \qquad \begin{aligned}
   C^\infty_{\text{even}}(S^3_1) &:= 
    \{  F\in C^\infty(S^3_1) \mid F=F\action \sigma \}, \\
   C^\infty_{\text{odd}}(S^3_1) &:= 
    \{  F\in C^\infty(S^3_1) \mid F=-F\action \sigma \}. 
 \end{aligned} $$
We call $h\in C^\infty_{\text{even}}(\bS^2)$ an {\it even function}, and so on. 
We define the maps 
$\pye : C^\infty(\bS^2)\to \R$ and 
$\xye : C^\infty(S^3_1)\to \R$ by 
$$ \begin{aligned}
 \pye(h)&= \int_{\bS^2}h \, \w_{\bS^2} \qquad \text{for}\quad 
     h\in C^\infty(\bS^2),\\ 
 \xye(F)&= \int_{v\in \bS^2} F(0,v) \, \w_{\bS^2} 
  \qquad \text{for}\quad F\in C^\infty(S^3_1) \end{aligned} $$
where $\w_{\bS^2}$ is the volume form on $\bS^2$. 
We set 
$$ \begin{aligned}
 C^\infty_*(\bS^2) &:= \{ h\in C^\infty(\bS^2) \mid \pye(h)=0 \}, \\
 C^\infty_{\text{even}\,*}(\bS^2) &:= 
 \{ h\in C^\infty_{\text{even}}(\bS^2) \mid \pye(h)=0 \}, 
 \end{aligned} \qquad \begin{aligned}
  C^\infty_*(S^3_1) &:=  \{ F\in C^\infty(S^3_1) \mid \xye(F)=0 \}, \\
  C^\infty_{\text{even}\,*}(S^3_1) &:= 
 \{ F\in C^\infty_{\text{even}}(S^3_1) \mid \xye(F)=0 \}. 
 \end{aligned} $$
Let us denote the space of real valued constant functions by $\R$.  
Then we obtain the natural decompositions 
$$ C^\infty(\bS^2)= \R \oplus C^\infty_{\text{even}\,*}(\bS^2) 
 \oplus C^\infty_{\text{odd}}(\bS^2) \quad \text{and} \quad
 C^\infty(S^3_1)= \R \oplus C^\infty_{\text{even}\,*}(S^3_1) 
 \oplus C^\infty_{\text{odd}}(S^3_1) $$
given by 
$ h = \pye(h) + \( \frac{1}{2}(h+h\action \a) - \pye(h)\) + 
 \frac{1}{2}(h-h\action \a) $ 
and so on.

\paragraph{Transforms}
 
We define linear transforms
$R, Q : C^\infty(\bS^2)\to C^\infty(S^3_1)$ by 
\begin{equation} \label{eq:basic_transform}
  R h (t,y) := \frac{1}{2\pi} \int_0^{2\pi} h(\g_{(t,y)}(\p)) d\p, 
\end{equation}
\begin{equation}
  Q h (t,y) := \frac{1}{2\pi} \int_{\Omega_{(t,y)}} h(u) \w_{\bS^2}, 
  \hspace{5.1mm}
\end{equation}
where $\g_{(t,y)}(\p)$ is given by \eqref{eq:parameter_of_circle}. 
Of course, $R$ is well-defined by \eqref{eq:basic_transform}
without depending on the choice of vectors $\{y^\perp_1,y^\perp_2 \}$. 
By definition we obtain 
$$ R (C^\infty(\bS^2)) \subset C^\infty_{\text{even}}(S^3_1), \qquad 
   Q (C^\infty_*(\bS^2))  \subset C^\infty_{\text{odd}}(S^3_1). $$

Restricting $R$ and $Q$ on the neck sphere $S^2_o\cong\bS^2$, 
we also define linear transforms 
$\cR, \cQ : C^\infty(\bS^2)\to C^\infty(\bS^2)$ by 
\begin{equation}
 \cR h (y):= Rh(0,y) \qquad \text{and} \qquad \cQ h (y):= Qh(0,y).  
\end{equation}
The transform $\cR$ is called the {\it Funk transform} (cf.\cite{bib:Funk}) 
or the {\it spherical Radon transform}. 
See \cite{bib:Guillemin76,bib:Helgason} for the detail of the 
(spherical) Radon transform and the related topics. 
The inverse problem for the (spherical) Radon transform is a classical problem, 
and there are a number of works on this subject. 
Recent development on the inverse problem concerning 
the Radon transform or related transforms are found in 
\cite{bib:Li-Song,bib:Rubin02} and the references in them. 
On the other hand, the transform $\cQ$ seems to be paid few attentions. 
We will study the inverse problem for the transform $\cQ$ 
in the next section and in Appendix A. 
We will apply the results of this study to solve the wave equation.

\begin{Lem} \label{Lem:del_R}
 For any smooth function $h$ on $\bS^2$, the following equation holds: 
 $$ \pd{t} Rh(t,y)= - Q\Delta_{\bS^2} h(t,y). $$  
\end{Lem}
\begin{proof}
Since we can vary $ t$ fixing the frame $\{y^\perp_1,y^\perp_2,y\}$ 
of $\R^3$, we obtain 
$$ \frac{\del}{\del t} (\g_{(t,y)}(\p))= - (\cosh t)^{-1} \nu(\p) $$ 
where $\nu(\p)$ is the unit normal vector field  along $\g_{(t,y)}(\p)$ 
directing outside of the domain $\Omega_{(t,y)}$.
Let $dm$ be the measure on $\del\Omega_{(t,y)}$ induced by the 
standard metric on $\bS^2$, then we have 
$dm=(\cosh t)^{-1}d\p$. Hence we obtain 
$$ \pd{t} R h (t,y)
 = - \frac{1}{2\pi}\int_{\del\Omega_{(t,y)}}
    (\nabla h)\cdot \nu \, dm 
 = - \frac{1}{2\pi}\int_{\Omega_{(t,y)}} (\Delta_{\bS^2}h)\, \w_{\bS^2} 
 = - Q\Delta_{\bS^2} h(t,y) $$
by the divergence formula. 
\end{proof}

\begin{Lem} \label{Lem:del_Q}
 For any smooth function $h$ on $\bS^2$, the following equation holds: 
 $$ \pd{t} Qh(t,y) = - (\cosh t)^{-2} Rh (t,y). $$  
\end{Lem}
\begin{proof}
 We fix $y\in\bS^2$ and take vectors $y^\perp_1,y^\perp_2\in \bS^2$ 
 so that $\{y^\perp_1,y^\perp_2,y\}$ gives an oriented orthonormal basis 
 on $\R^3$. We use a spherical coordinate $(\c,\p) \in (0,\pi)\times(0,2\pi)$ 
 on $\bS^2$ defined by 
 $$ (\c,\p) \longmapsto 
  u(\c,\p)= 
  \sin\c\cos\p \ y^\perp_1 + \sin\c\sin\p \ y^\perp_2 + \cos\c \ y. $$ 
 Then we have 
 $ \Omega_{(t,y)}= \{ (\c,\p) \in \bS^2 \mid 0\le \c \le \a \} $ 
 where $\a$ is the real variable defined by $\cos\a=\tanh t$. 
 In this coordinate, noticing $\w_{\bS^2}=\sin\c\, d\c\wedge d\p$, 
 $$ Qh(t,y)= \int_{\Omega_{(t,y)}} h(u(\c,\p)) \sin\c\, d\c\wedge d\p 
   = \int_0^{2\pi} \left[ \int_0^\a h(u(\c,\p)) \sin\c\, d\c \right] d\p. $$
 Since $\pd{t}= - (\cosh t)^{-1}\pd{\a}$, we obtain 
 $$ \pd{t} Qh(t,y)= \frac{-1}{\cosh t} \int_0^{2\pi} h(u(\a,\p)) \sin\a\, d\p 
   = \frac{-1}{\cosh^2 t} \int_0^{2\pi} h(u(\a,\p)) d\p
   = - \frac{Rh(t,y)}{\cosh^2 t} $$
 as required. 
\end{proof}

\begin{Rem}
 We can check $Q(1)(t,y)={\rm Area}(\Omega_{(t,y)})=1-\tanh t$ 
 by using the above coordinate $(\c,\p)$.  
\end{Rem}

\begin{Prop} \label{Prop:Rh_solves_MPPE}
 For any smooth function $h$ on $\bS^2$, the induced function $f:=R h$ 
 on $S^3_1$ solves the following hyperbolic partial differential equation: 
 \begin{equation} \label{eq:MPPE}
  Lf := \( -\frac{\del^2}{\del t^2} + (\cosh t)^{-2}\Delta_{\bS^2} \) f =0. 
 \end{equation}
\end{Prop}
\begin{proof}
First we claim that $R$ commutes with $\Delta_{\bS^2}$. 
Actually, if we fix $t\in \R$, the transform $M^t: h \longmapsto Rh(t,\cdot)$ is 
$SO(3)$-equivariant, hence $M^t$ commutes with $\Delta_{\bS^2}$ 
by Theorem \ref{Thm:decomp.of_L2S2} in Appendix A. 
Thus $R$ commutes with $\Delta_{\bS^2}$. 
Then, by above Lemmas, 
$$ \frac{\del^2}{\del t^2} R h (t,y) 
 = - \pd{t} Q\Delta_{\bS^2} h(t,y) 
 = (\cosh t)^{-2} R \Delta_{\bS^2} h \, (t,y)
 = (\cosh t)^{-2} \Delta_{\bS^2} R h \, (t,y) $$
for any smooth function $h$ on $\bS^2$. 
Hence $f:=Rh$ solves \eqref{eq:MPPE}. 
\end{proof}

\begin{Lem} \label{Lem:MPPE=>wave}
 Let $f$ be a smooth function on $S^3_1$ satisfying the equation 
 $Lf=0$. 
 If we put $V:=f_t$, then $V$ satisfies the wave equation $\Box V=0$. 
\end{Lem}
\begin{proof}
 Applying $\pd{t}$ on the equation $(\cosh t)^2 Lf=0$, we obtain 
 the equation \eqref{eq:WE}. 
\end{proof}

\begin{Prop} \label{Prop:Qh_solves_wave_eq}
 Let $h$ be a smooth function on $\bS^2$ satisfying $\pye(h)=0$. 
 Then $V :=Qh$ solves the wave equation $\Box V=0$. 
\end{Prop}
\begin{proof}
Since $\pye(h)=0$, there exists a smooth function $\tilde h$
on $\bS^2$ satisfying $h=-\Delta_{\bS^2}\tilde h$. 
If we put $f:=R\tilde h$, then $V=-Q\Delta_{\bS^2}\tilde h=f_t$ 
by Lemma \ref{Lem:del_R}. 
On the other hand, $Lf=0$ by Proposition \ref{Prop:Rh_solves_MPPE}, 
so $\Box V=\Box f_t=0$ by  Lemma \ref{Lem:MPPE=>wave}. 
 \end{proof}

\begin{Rem}
We call a function $h\in C^\infty_*(\bS^2)$ a {\it generating function} 
in the sense that $h$ induces a solution of $\Box V=0$ or $Lf=0$. 
\end{Rem}

\section{Oddness and the inverse problem}
\label{Sect:Inverse}

In this section, we investigate the inverse problem for the transform 
$R$ and $Q$. The goal is the following.  

\begin{Thm} \label{Thm:Inverse}
 \begin{enumerate}
 \item 
  Let $V$ be a smooth function on $S^3_1$ which solves the wave equation 
  $\Box V=0$.   Suppose that $V\to 0$ and $V_t\to 0$ as $t\to\pm\infty$ 
  uniformly for $y\in \bS^2$. 
  Then $V$ is odd, and there exists a unique smooth function 
  $h\in C^\infty_*(\bS^2)$ satisfying $V=Qh$. 
 \item 
  Let $f$ be a smooth function on $S^3_1$ which solves the equation $Lf=0$. 
  Suppose that 
  there exist $h_\pm(y)\in C^\infty(\bS^2)$ and that 
  $f(t,y)\to h_\pm(y)$ and $f_t, f_{tt}\to 0$ as $t\to \pm\infty$ 
  uniformly for $y\in\bS^2$. 
  Then $f$ is even, and $f=Rh_+$ holds.  
  Moreover, if $f\in\C^\infty_*(S^3_1)$ then $h_\pm \in C^\infty_*(\bS^2)$. 
 \end{enumerate}
\end{Thm}

\paragraph{Inverse problem for $\cR$ and $\cQ$}

First, we study the transforms $\cR$ and $\cQ$. 
By definition, $\cR$ and $\cQ$ are the identities on the constant functions 
$\R\subset C^\infty(\bS^2)$. We also have 
$\cR(C^\infty_{\text{odd}}(\bS^2))=
 \cQ(C^\infty_{\text{even}\, *}(\bS^2))=0$. 
The following bijectivity is, however, rather non-trivial. 

\begin{Prop} \label{Prop:bijectivity_of_the_transforms}
 Both of the following transforms are bijective: 
 \begin{enumerate}
  \item $\cR : C^\infty_{\text{even}\, *}(\bS^2) 
    \longrightarrow C^\infty_{\text{even}\, *}(\bS^2)$, 
  \item $\cQ : C^\infty_{\text{odd}}(\bS^2) 
    \longrightarrow C^\infty_{\text{odd}}(\bS^2)$. 
 \end{enumerate}
 Hence we obtain 
 $$ \ker\{\cR:C^\infty(\bS^2)\to C^\infty(\bS^2)\}
  = C^\infty_{\text{odd}}(\bS^2), \qquad 
   \ker\{\cQ:C^\infty(\bS^2)\to C^\infty(\bS^2)\}
  = C^\infty_{\text{even}\, *}(\bS^2). $$
\end{Prop}

The above bijectivity of $\cR$ on $C^\infty_{\text{even}\, *}(\bS^2)$
is first noticed by P.~Funk \cite{bib:Funk}. 
There is an explicit inversion formula of $\cR$ or its generalization, 
which we can find in the textbook by S.~Helgason \cite{bib:Helgason}. 
On the other hand, just for the purpose to verify the bijectivity of $\cR$, 
V.~Guillemin's method is reasonable (Appendix A of \cite{bib:Guillemin76}). 
We give the proof of the bijectivity of $\cQ$ on 
$C^\infty_{\text{odd}}(\bS^2)$ in the Appendix A 
by a similar argument as the Guillemin's.

\paragraph{Key lemma} 
The key to prove Theorem \ref{Thm:Inverse} is to verify the oddness and evenness 
for the initial values $V|_{t=0}$ and $V_t|_{t=0}$, which will be shown in 
Lemma \ref{Lem:oddness_of_V}. Before this, we first notice the following. 

\begin{Lem} \label{Lem:V_in_C_*}
 For each function $V\in C^\infty(S^3_1)$, 
 we define a function $I(\t)$ on $\t\in \R$ by  
 \begin{equation} \label{eq:I(tau)}
 I(\t):=\frac{\cosh^2\! \t}{2\pi} \int_{\bS^2} (V_t|_{t=\t})\, \w_{\bS^2}. 
 \end{equation}
 If \, $V$ \! solves the wave equation $\Box V=0$, 
 then $I(\t)$ is independent with $\t\in\R$. 
\end{Lem}
\begin{proof}
 Let $p: S^3_1=\R\times\bS^2 \to \R$ be the projection, 
 and we notice to the interval $(t_1,t_2)\subset \R$. 
 If $V$ satisfies $\Box V= *\,d*d\, V=0$, then we obtain 
 $$ 0= \frac{1}{2\pi}\int_{p^{-1}(t_1,t_2)} d*d\, V 
 = \frac{-1}{2\pi} \int_{p^{-1}(t_2)-p^{-1}(t_1)} V_t \cosh^2\! t \ \w_{\bS^2} 
 = - I(t_2) + I(t_1) $$ 
 Hence $I(\t)$ does not depend on $\t\in\R$. 
\end{proof}

\begin{Lem} \label{Lem:oddness_of_V}
 Let $V$ be a smooth function on $S^3_1$ which solves 
 the wave equation $\Box V=0$. 
 Suppose that $V\to 0$ and $V_t\to 0$ as $t\to\pm\infty$ 
 uniformly for $y\in \bS^2$. 
 Let $\psi(y):=V(0,y)$ and $\x(y):=V_t(0,y)$. 
 Then $\psi\in C^\infty_{\text{odd}}(\bS^2)$ and 
 $\x \in C^\infty_{\text{even}\, *}(\bS^2)$.   
\end{Lem}
\begin{proof}
 We fix $y\in\bS^2$ and take a coordinate $(\c,\p)$ on $\bS^2$ 
 similarly as in the proof of Lemma \ref{Lem:del_Q}. 
 We use the coordinate $(t,\c,\p) \in \R\times\bS^2\cong S^3_1$, 
 and we put 
 $$ \begin{aligned}
 \Omega'_{(t,y)} &:= \{t\}\times \Omega_{(t,v)} \subset S^3_1, \\
 M_y(t_1,t_2) &:= \{ (t,\c,\p) \in S^3_1 
   \mid t_1<t<t_2, \cos\c>\tanh t \}
   = \cup_{t\in (t_1,t_2)} \Omega'_{(t,y)}, \\
 \Sigma_y(t_1,t_2) & := \{ (t,\c,\p) \in S^3_1 
   \mid t_1<t<t_2, \cos\c=\tanh t \}
   = \cup_{t\in (t_1,t_2)} \del \Omega'_{(t,y)}.
  \end{aligned} $$  
 Notice that 
 $$ \begin{aligned} 
  \del M_y(t_1,t_2) &= \Sigma_y(t_1,t_2) \cup \Omega'_{(t_2,y)} 
    \cup (-\Omega'_{(t_1,y)}) \\ 
  \del\Sigma_y(t_1,t_2) &= \del \Omega'_{(t_2,y)}\cup 
    (-\del\Omega'_{(t_1,y)}). \end{aligned} $$ 

 Now let $V$ be a function on $S^3_1$ 
 as in the statement and $\t$ be a positive real variable. 
 Since $\Box V=*\, d*dV=0$, 
 integrating on $M_y(t_1,t_2)$, we obtain 
 \begin{equation} \label{eq:total_integral}
   0 = \int_{M_y(t_1,t_2)} d*d\, V  
      = \int_{\Sigma_y(t_1,t_2)+ \Omega'_{(t_2,y)} - \Omega'_{(t_1,y)}} *\, d\, V. 
 \end{equation}
 For $i=1,2$, let $\a_i\in (0,\pi)$ be the real variable 
 defined by $\cos \a_i = \tanh t_i$. 
 To calculate the integral over $\Sigma_y(t_1,t_2)$, we introduce a 
 real coordinate $(a,b)\in (\a_2,\a_1)\times (0,2\pi)$ on $\Sigma_y(t_1,t_2)$
 by the embedding  
 $ j :  (\a_2,\a_1)\times (0,2\pi) \rightarrow \Sigma_y(t_1,t_2) $ 
 defined by $(t,\c,\phi)=(t(a),a,b)$ where $\cos a = \tanh t(a)$. 
 Then we obtain 
 $$ \begin{aligned} 
  *\, d\, V &=-V_t \cosh^2\! t \ \w_{\bS^2}- dt \wedge \check *\, \check d\, V 
   = \( -(\sin a)^{-1} V_t + V_\c \) da\wedge db 
   = - \PD{V}{a} da\wedge db \\
   &= - d(Vdb). \end{aligned} $$ 
 Hence 
 $$ \begin{aligned} 
  \frac{1}{2\pi} \int_{\Sigma_y(t_1,t_2)} *\, d\, V
  &= \frac{-1}{2\pi}  \int_{\Sigma_y(t_1,t_2)} d(Vdb) 
  = \frac{-1}{2\pi} \int_{\del \Omega'_{(t_2,y)}-\del \Omega'_{(t_1,y)}}  V d\phi 
   \\[1eX] 
  & = -R(V|_{t=t_2}) (t_2,y) + R(V|_{t=t_1}) (t_1,y) .  
 \end{aligned} $$ 
 On the other hand, for each $\t\in\R$ we have 
 $$ \frac{1}{2\pi} \int_{\Omega'_{(\t,y)}} *\, d\, V 
  = \frac{-1}{2\pi} \int_{\Omega'_{(\t,y)}} V_t \cosh^2\! t \ \w_{\bS^2} 
  = -\cosh^2 \!\t \ Q(V_t|_{t=\t}) (\t,y) $$ 
 Hence by \eqref{eq:total_integral}, we see that the quantity 
 \begin{equation} \label{eq:energy}
  E(y) := R(V|_{t=\t}) (\t,y) + \cosh^2\!\t \ Q(V_t|_{t=\t}) (\t,y) 
 \end{equation}
 does not depend on $\t\in \R$. 

 Now we claim $E(y)\equiv 0$. Notice that 
 \begin{equation} \label{eq:Qxi}
  \cosh^2\! \t \left| \int_{\Omega'_{(\t,y)}} V_t\, \w_{\bS^2} \right| 
	 \le \cosh^2\!\t \, {\rm Area}(\Omega_{(\t,y)}) \cdot \!
	  \max_{\, u\in\Omega'_{(\t,y)}} |V_t(u)| 
	  \ \le \max_{\, u\in\Omega'_{(\t,y)}} |V_t(u)| 
 \end{equation}
 since ${\rm Area}(\Omega_{(\t,y)})=1-\tanh \t$. 
 Hence we obtain $\lim_{\t\to +\infty} [\cosh^2\!\t \, Q(V_t|_{t=\t}) (\t,y)]=0$. 
 by the convergence of $V_t$. 
 On the other hand, by the convergence of $V$, we also have  
 $\lim_{\t\to\pm\infty} [R(V|_{t=\t}) (\t,y)] =0$. 
 Thus, by taking the limit $\t\to +\infty$ on \eqref{eq:energy}, 
 we obtain $E(y)\equiv 0$ as required. 
 
 Next notice that 
 \begin{equation} \label{eq:decomposing_I}
   I = \cosh^2\!\t \ Q(V_t|_{t=\t}) (\t,y) + 
  \cosh^2\!\t \ Q(V_t|_{t=\t}) (-\t,-y) 
 \end{equation} 
 where $I$ is the quantity defined in \eqref{eq:I(tau)}.
 If we take the limit $\t\to -\infty$, then the second term of the right hand side of 
 \eqref{eq:decomposing_I} vanishes by the similar argument as above. 
 Hence we obtain 
 $$ I=\lim_{\t\to -\infty} \left[ \cosh^2\!\t \ Q(V_t|_{t=\t}) (\t,y) \right]. $$
 Thus, by taking the limit $\t\to -\infty$ on \eqref{eq:energy}, we obtain 
 $ I=0$. This means $V_t|_{t=\t}\in C^\infty_*(\bS^2)$ for any $\t\in \R$. 
 
 Finally, evaluating $\t=0$ to  \eqref{eq:energy}, 
 we obtain $\cR \psi +\cQ \x=0$. 
 Recall that $\cR \psi \in C^\infty_{\text{even}}(\bS^2)$ and 
 $\cQ \x \in \R\oplus C^\infty_{\text{odd}}(\bS^2)$ 
 by Proposition \ref{Prop:bijectivity_of_the_transforms}. 
 Further, since $\x=V_t|_{t=0}\in C^\infty_*(\bS^2)$ by the above argument, 
 we have $\cQ \x \in C^\infty_{\text{odd}}(\bS^2)$. 
 Thus we obtain  $\cR \psi =\cQ \x =0$. 
 Hence $\psi \in C^\infty_{\text{odd}}(\bS^2)$ and  
 $\x \in C^\infty_{\text{even}\, *}(\bS^2)$
 by Proposition \ref{Prop:bijectivity_of_the_transforms}.  
\end{proof}

\paragraph{Proof of Theorem \ref{Thm:Inverse}} 

Let $V$ be as in the statement {\it 1}, and let 
$\psi(y):=V(0,y)$ and $\x(y) := V_t(0,y)$. 
By Lemma \ref{Lem:oddness_of_V}, we have 
$\psi\in C^\infty_{\text{odd}}(\bS^2)$ and 
$\x\in C^\infty_{\text{even}\, *}(\bS^2)$. 
Then by Proposition \ref{Prop:bijectivity_of_the_transforms}, 
there exist smooth functions 
$h_{\rm odd}\in C^\infty_{\rm{odd}}(\bS^2)$ and 
$h_{\rm even}\in C^\infty_{\rm{even}\, *}(\bS^2)$ satisfying 
\begin{equation} \label{eq:inverse_of_V}
 \psi = \cQ h_{\rm odd}, \qquad 
 \x = - \cR h_{\rm even}. 
\end{equation}
Now let us put $h:=h_{\rm even}+h_{\rm odd}$ and $\tilde V:= Qh$. 
Since $h\in C^\infty_*(\bS^2)$, $\tilde V$ is a solution of $\Box \tilde V=0$ 
by Proposition \ref{Prop:Qh_solves_wave_eq}. Moreover by construction 
$$ \tilde V(0,y)= \cQ h (y) = \psi(y), \qquad 
  \tilde V_t(0,y)= \left[ \pd{t} Qh(t,y) \right]_{t=0} = - \cR h(y) = \x(y). $$ 
Hence $V$ and $\tilde V$ satisfies the same initial condition, so by 
the uniqueness theorem for the initial value problem of 
hyperbolic partial differential equations (see \cite{bib:CH}), we obtain 
$V=\tilde V$. Hence $V=Qh$ and $V$ turns out to be odd. 
The uniqueness of $h$ is obvious by the relation \eqref{eq:inverse_of_V}. 

Next let $f$ be as in the statement {\it 2}. If we put $V:=f_t$, then 
$V$ satisfies the conditions in the statement {\it 1}. 
Hence $V$ is odd. 
If we decompose $f$ as $f= f_{\text{even}} + f_{\text{odd}}$ so that 
$f_{\text{even}}\in C^\infty_{\text{even}}(S^3_1)$ and  
$f_{\text{odd}}\in C^\infty_{\text{odd}}(S^3_1)$, 
then $V_t=f_t=(f_{\text{odd}})_t + (f_{\text{even}})_t$ gives the 
decomposition of $V$ satisfying 
$(f_{\text{odd}})_t\in C^\infty_{\text{even}}(S^3_1)$ and  
$(f_{\text{even}})_t\in C^\infty_{\text{odd}}(S^3_1)$. 
Since $V$ is odd, we obtain $(f_{\text{odd}})_t=0$. 
Hence $f_{\text{odd}}=0$ and $f$ is even. 

Let us put $\vp(y):= f(0,y)$ and $\psi(y):= f_t(0,y)$. 
Then similar as the above argument, 
there is a unique smooth function $\tilde h$ on $\bS^2$ which satisfies 
$$ \vp(y)= \cR \tilde h (y), \qquad \psi(y)=-\cQ \Delta_{\bS^2}\tilde h (y). $$
For this function $\tilde{h}$, we obtain $f=R \tilde h$. 
By definition of $R$, 
$$h_+(y)=\lim_{t\to\infty}f(t,y)=\lim_{t\to\infty}R\tilde h (t,y)=\tilde h(y).$$ 
Hence $f=R h_+$ as required. 
If $f\in C^\infty_*(S^3_1)$, then $\vp\in C^\infty_{\text{even}\, *}(\bS^2)$ 
and we obtain $h_+=\tilde h\in C^\infty_*(\bS^2)$ by the construction.  
Since $f$ is even, $h_-(y)=h_+(-y)\in C^\infty_*(\bS^2)$. 
\qed

\paragraph{Tameness at infinity}
By Theorem \ref{Thm:Inverse} and its proof, we can paraphrase the condition of 
the `tameness at infinity' for $V$ in the following way. 

\begin{Cor} 
 Let $V\in C^\infty(S^3_1)$ be a solution of the wave equation $\Box V=0$. 
 then the following conditions are equivalent: 
 \begin{enumerate}
  \item  $V(t,y)\to 0$ and $V_t(t,y)\to 0$ as $t\to\pm\infty$ 
   uniformly for $y\in \bS^2$, 
   \item $V$ is odd and $I=0$, and  
   \item $\psi(y):=V(0,y) \in C^\infty_{\text{odd}}(\bS^2)$ and 
      $\xi(y):=V_t(0,y) \in C^\infty_{\text{even}\, *}(\bS^2)$. 
 \end{enumerate}
\end{Cor}
\begin{proof}
 The statement {\it 1} $\Rightarrow$ {\it 2} follows 
 from Theorem \ref{Thm:Inverse},  
 and {\it 2} $\Rightarrow$ {\it 3} is obvious.
 Now let us assume {\it 3}. 
 As in the proof of Theorem \ref{Thm:Inverse}, 
 we have $V=Qh$ for $h=h_{\text{even}}+h_{\text{odd}}$ where 
 $h_{\text{even}}\in C^\infty_{\text{even}\, *}(\bS^2)$ and 
 $h_{\text{odd}}\in C^\infty_{\text{odd}}(\bS^2)$ are defined by 
 \eqref{eq:inverse_of_V}. 
 Then we can check that $V=Qh$ and $V_t=-(\cosh t)^{-2} Rh$ 
 uniformly converge to zero 
 as $t\to\pm\infty$. Thus {\it 3} $\Rightarrow$ {\it 1} holds. 
\end{proof}

Similarly, we obtain the following corollary of which the proof is omitted. 

\begin{Cor}
 Let $f\in C^\infty(S^3_1)$ be a solution of the equation $Lf=0$. 
 then the following conditions are equivalent: 
 \begin{enumerate}
  \item There exist smooth functions $h_{\pm}(y)\in C^\infty(\bS^2)$ 
    such that $f(t,y)\to h_\pm(y)$ and $f_t,f_{tt}\to 0$ as $t\to \pm\infty$ 
	uniformly for $y\in\bS^2$, 
  \item $f$ is even, and 
  \item $\varphi(y):=f(0,y)$ is even and $\psi(y):=f_t(0,y)$ is odd. 
 \end{enumerate}
\end{Cor}

\paragraph{Rigidity theorem}

Let $S^3_1/\Z_2$ be the quotient space of $S^3_1$ by the involution $\s$. 
Notice that $S^3_1/\Z_2$ is not space-time-orientable. 
Since the operator $\Box$ on $S^3_1$ is $\sigma$-invariant, we can 
define the wave equation $\Box V=0$ on $S^3_1/\Z_2$. 
Now let us use the coordinate $\{(t,y)\in \R\times \bS^2 \mid t>0\}$ on the 
open set $\{[t,y]\in S^3_1/\Z_2 \mid t\neq 0\}$. 
Then, as a trivial consequence of Theorem \ref{Thm:Inverse}, 
we obtain the following rigidity theorem. 

\begin{Cor} \label{Cor:rigidity_of_waves}
 Let $V$ be a solution of the wave equation $\Box V=0$ on $S^3_1/\Z_2$. 
 Suppose $V,V_t\to 0$ as $t\to\infty$ uniformly for $y$, then $V\equiv 0$. 
\end{Cor}

We remark that this type of rigidity theorem is also 
found in \cite{bib:LM02} or \cite{bib:LM05}. 
For example in \cite{bib:LM05}, 
it is shown that the standard self-dual indefinite metric 
on the non-space-time-orientable space $(S^2\times S^2)/\Z_2$ 
is rigid in the space of self-dual metrics.

\section{Monopole equation}
\label{Sect:Monopole_equation}

In this section, we investigate the monopole equation 
over the de Sitter space $S^3_1$. 
We show that any gauge equivalence class of monopole solutions is 
obtained from a solution of $Lf=0$ which we call the {\it monopole potential}. 
Then, applying Theorem \ref{Thm:Inverse}, 
we establish a one-to-one correspondence between 
generating functions $h\in C^\infty_*(\bS^2)$ and gauge equivalent 
classes of monopole solutions on $S^3_1$ which are tame at infinity. 
Further, we introduce the notion of 
{\it admissible monopoles} by which we can construct $S^1$-invariant 
self-dual metrics on $S^2\times S^2$.

\paragraph{Tod-Kamada ansatz} 
Here we review the construction of self-dual metrics on $S^2\times S^2$ 
given by Tod or Kamada, 
following Kamada's formulation. 

The basic construction is the following. 

\begin{Prop}[Kamada\cite{bib:Kamada05}] \label{Prop:Kamada_ansatz}
Let $V$ be a smooth positive function on $S^3_1$ such that $*dV/2\pi$ 
is a closed two-form on $S^3_1$ determining an integral class in $H^2(S^3_1;\R)$. 
Let $\cM\to S^3_1$ denote an $S^1$-bundle with connection one-form 
$\Theta$ with curvature form given by 
\begin{equation} \label{eq:Kamada_MPE} d\Theta= *dV. \end{equation}
Then $g_{V,\Theta}:= -V^{-1}\Theta\otimes\Theta + V g_{S^3_1}$ 
is a self-dual metric on $\cM$ of signature $(--++)$ 
with respect to a suitable orientation on $\cM$. 
\end{Prop}

Now we study the case when $*dV$ is exact, i.e. 
when the $S^1$-bundle $\cM\to S^3_1$ is trivial.  
In this case, we write as 
$\cM\simeq S^1\times S^3_1=\{(s,t,y)\in S^1\times\R\times \bS^2\}$ 
where $s$ is the fiber coordinate and $S^3_1=\{ (t,y)\in \R\times \bS^2\}$. 
The total space $\cM$ is naturally compactified to 
$\bar\cM:=S^2\times S^2$ by the embedding 
$\cM\hookrightarrow S^2\times S^2 : (s,t,y)\mapsto (x,y)$ where 
\begin{equation} \label{eq:x<->(s,t)}
 x^1=\frac{\cos s}{\cosh t}, \quad 
 x^2=\frac{\sin s}{\cosh t}, \quad 
 x^3=\tanh t. 
\end{equation}
In other words, $\cM$ is obtained as the free part of the $S^1$-action on 
$S^2\times S^2$ defined by 
\begin{equation} \label{eq:S^1-action}
\a\cdot(x,y)=(R(\a)x,y), \quad \text{where} \ \ \a\in S^1
 \ \ \text{and} \ \ 
 R(\a)= \begin{pmatrix} 
  \cos\a& -\sin\a &0 \\ \sin\a &\cos\a &0 \\ 0&0&1\end{pmatrix}. 
\end{equation}
If we put $\e:=(0,0,1)\in S^2$ and 
$S_\pm:=\{\pm\e\}\times S^2 \subset S^2\times S^2$, 
then the disjoint union $S_+\sqcup S_-$ coincides to 
the fixed point set of the above $S^1$-action, and we have 
$\cM=(S^2\times S^2) \backslash(S_+\sqcup S_-)$. 
Let us introduce variables $r:=e^t$ and $q:=e^{-t}$, 
then $(s,r)$ and $(s,q)$ give the polar coordinates 
on the open neighborhoods of $-\e\in S^2$ and $\e\in S^2$ respectively. 

\begin{Prop}[Kamada\cite{bib:Kamada05}] \label{Prop:Kamada_compactification}
 Let $(V,\Theta)$ be a smooth solution of \eqref{eq:Kamada_MPE} 
 such that $V>0$ and $*dV$ is an exact two-form. Then the metric 
 $\bar{g}_{V,\Theta}:= (\cosh t)^{-2} g_{V,\Theta}$ 
 on $\cM$ extends smoothly to 
 the compactification $\bar\cM = S^2\times S^2$ 
 if and only if there exist smooth functions  $F_ +$ and $F_-$ 
 on $\R\times \bS^2$ in variables $r^2,q^2$ and $y$ such that 
 \begin{equation} \label{eq:Kamada_compactification}
  V=1+r^2 F_-(r^2,y) \quad \text{and} \quad
  V=1+q^2 F_+(q^2,y),
 \end{equation}
 as $r\to +0$ and as $q\to +0$ respectively. 
\end{Prop}

If $*dV$ is exact, 
$\Theta$ is written as $\Theta=ds+A$ using a one-form $A$ on $S^3_1$. 
Then the equation \eqref{eq:Kamada_MPE} is written as 
\begin{equation} \label{eq:MPE}
dA=* dV 
\end{equation}
which we call the {\it monopole equation}. 
We call a solution $(V,A)$ of \eqref{eq:MPE} a 
{\it monopole solution} or simply a {\it monopole}. 
We write as $\bar{g}_{V,A}=\bar{g}_{V,\Theta}$ where 
$\Theta=ds+A$, and we also use the notation $\bar{g}_{V,A}$ for its 
compactification. 
Notice that if $(V,A)$ is a monopole then $V$ satisfies the 
wave equation $\Box V= *d*d V=0$. 

The simplest solution of the monopole equation satisfying the condition 
\eqref{eq:Kamada_compactification} is given by 
$(V,A)=(1,0)$, which we call the {\it trivial monopole}. In this case, 
the self-dual indefinite metric induced on $S^2\times S^2$ is 
the standard indefinite metric, i.e. the product metric $g_0=\pi_1^*h-\pi_2^*h$ 
where $\pi_i:S^2\times S^2\to S^2$ is the $i$-th projection and $h$ is the 
standard metric on $S^2$.

Tod's or Kamada's examples of self-dual indefinite metrics are obtained by 
constructing explicit solutions of \eqref{eq:MPE}. We deal with these examples 
in the last part of this section.

\paragraph{Monopole potential}
Now we show that any monopole solution is essentially arisen from a function 
$f\in C^\infty_*(S^3_1)$ satisfying $Lf=0$ where 
$L$ is the partial differential operator defined in \eqref{eq:MPPE}. 
We call such $f$ the {\it monopole potential}.

For each real valued function $\p\in C^\infty(S^3_1)$, 
the transform of monopoles 
$$ (V,A) \longmapsto (V,A + d\p) $$
is called the {\it gauge transform}. 
Notice that $\Phi^* \bar{g}_{V,A}=\bar{g}_{V,A+d\p}$ where 
$\Phi=e^{i\p} : \cM \to \cM$ is the gauge transform on the 
$S^1$-bundle $\cM\to S^3_1$. 

\begin{Prop} \label{Prop:gauge_fixing}
 Let $(V,A)$ be any monopole on $S^3_1$. 
 Then, changing $(V,A)$ by a gauge transform, we can assume 
 $({\rm 1}^\circ)$ $A(\del_t)=0$ and 
 $({\rm 2}^\circ)$ $\check d \, \check * \, A=0$, 
 where $\check d$ is the fiberwise exterior derivative 
 and $\check *$ is the fiberwise Hodge's operator 
 on the $\bS^2$-bundle $S^3_1\to\R$.  
 Furthermore such $(V,A)$ is unique in the gauge equivalence class. 
\end{Prop}
\begin{proof}
 Let $(V,A)$ be any monopole on $S^3_1$. 
 Let us write as $A=A_t dt +A_1$ so that $A_1$ is a 1-form without $dt$-part. 
 If we put $\p_1=-\int_0^t A_tdt$, then the one-form $A+ d\p_1$ does not 
 have $dt$-part. Hence we can assume that $A$ satisfies 
 the condition $1^\circ$ from the beginning. 
 
 By the monopole equation \eqref{eq:MPE}, $A=A_1$ satisfies $d*d\, A=0$. 
 Since $dA = dt \wedge \PD{A}{t} + \check d \, A, $
 we can write as 
 $$ 0 = d*d\, A \equiv d*\( dt\wedge \PD{A}{t} \) 
  \equiv \pd{t}\, (d*(dt\wedge A) )  
  \equiv -\pd{t}\, (\check d \, \check * \, A) 
  \quad \mod \ dt, $$ 
 where we applied the relation $*(dt\wedge \y)=- \check *\, \y$ which holds 
 for any 1-form $\y$ on $S^3_1$ without $dt$-part. 
 Thus the function 
 $\check d^{\, *} A := \check * \,\check d\,\check * A$
 does not depend on $t$. 
 
 Considering $\check d^{\, *} A$ as a function on $\bS^2$,
 we can take a smooth function $\check \p$ on $\bS^2$ satisfying 
 $\Delta_{\bS^2}\check \p =-\check d^{\, *} A$ 
 since $(\check d^{\, *} A) \, \w_{\bS^2} =
  \check * \, (\check d^{\, *} A) = 
  \check d \, \check *\, A$ is exact. 
 If we define a smooth function $\p$ as the pull back of $\check \p$ by the 
 projection $S^3_1\cong \bS^2\times\R \to \bS^2$, then we obtain 
 $ \check d\, \check * (A+ d\p) =0$. Hence 
 $A'=A+ d\p$ satisfies the conditions $1^\circ$ and $2^\circ$. 

 Now we prove the uniqueness. 
 Suppose that there is a function $\p$ on $S^3_1$ such that 
 both $(V,A)$ and $(V,A+d\p)$ are the monopoles 
 satisfying $1^\circ$ and $2^\circ$. 
 Then the monopole $(0,d\p)$ also satisfies $1^\circ$ and $2^\circ$. 
 By condition $1^\circ$, $\p$ is independent of $t$. 
 Hence $d\p=\check d \p$. 
 Together with the condition $2^\circ$, we obtain $d\p=0$. 
 So the uniqueness follows. 
\end{proof}

\begin{Prop} \label{Prop:monopole_potential}
 Let $(V,A)$ be a monopole on $S^3_1$. 
 Suppose that $A$ satisfies 
 $({\rm 1}^\circ)$ $A(\del_t)=0$ and 
 $({\rm 2}^\circ)$ $\check d \, \check * \, A=0$.  
 Then there exists a unique function $f\in C^\infty_*(S^3_1)$ 
 satisfying (i) $V=\del_t f$ and (ii) $A=-\check * \, \check d \, f$. 
 Moreover $f$ satisfies the equation $Lf=0$, where $L$ is the partial differential 
 operator defined in \eqref{eq:MPPE}. 
\end{Prop}
\begin{proof}
 Let $(V,A)$ be a monopole on $S^3_1$ satisfying ${\rm 1}^\circ$ 
 and ${\rm 2}^\circ$. 
 We first claim that there is a smooth function $F$ on $S^3_1$ such that 
 $A=-\check *\, \check d \, F$. 
 Such a function is obtained, for example, 
 by putting $F(t,y):= \int_{(t,o)}^{(t,y)} \check *\, A$ where 
 $o\in \bS^2$ is a fixed point and the integral path is taken on the sphere
 $\{t\}\times \bS^2\subset S^3_1$. 
 Since $\bS^2$ is simply connected, and by the 
 condition ${\rm 2}^\circ$, 
 $F(t,y)$ is a well-defined smooth function. 
 By construction, the condition $A=-\check * \, \check d \, F$ holds. 

 Next we claim that $\check d \, (V-\del_t F)=0$. Actually, 
 $$ \begin{aligned} 
  \check d\, (V-\del_t F) &= dV - V_t dt - \check d\, F_t, \\
  * \, \check d\, (V-\del_t F) &= *dV + V_t \w_{\bS^2} + 
    dt \wedge \check *\, \check d\, F_t
    = \check d\, A + V_t \w_{\bS^2},
     \end{aligned} $$ 
 and $\check d\, A + V_t \w_{\bS^2}=0$ by the monopole equation. 
 Hence $\check d \, (V-\del_t F)=0$ as required, and this means 
 that $G(t):=V(t,y)-\del_t F(t,y)$ does not dependent on $y\in \bS^2$. 
 Thus, if we put $f(t,y)=F(t,y)+\int_0^t G(t)dt$, the condition (i) and (ii) are 
 satisfied. 
 The uniqueness of $f$ is obvious since the conditions (i) and (ii) characterize 
 $f$ up to constant. 

 The rest of the statement is directly follows from the monopole equation. Indeed, 
 $$ \begin{aligned}
   *\, dV &= *\, (f_{tt}\, dt + \check d\, f_t) 
       = - f_{tt} \cosh^2\!t \, \w_{\bS^2} 
	      - dt \wedge \check *\, \check d\, f_t, \\ 
   dA &= - d\, \check *\, \check d\, f 
     = - dt \wedge (\check *\, \check d\, f)_t - \check d\,\check *\, \check d\, f
	 = - dt \wedge \check *\, \check d\, f_t - (\Delta_{\bS^2}f)\, \w_{\bS^2}, 
   \end{aligned} $$
 hence $0=*\, dV - dA = (Lf)\cosh^2\!t \, \w_{\bS^2}$. 
\end{proof}

For monopoles which are tame at infinity, 
we obtain the following correspondence.

\begin{Thm} \label{Thm:3objects}
 There is a natural one-to-one correspondence between 
 the following objects:
 \begin{itemize}
  \item[1.] {\bf [generating functions]} \\
   smooth functions $h\in C^\infty_*(\bS^2)$, 
  \item[2.] {\bf [monopole potentials]} \\
   smooth functions $f\in C^\infty_*(S^3_1)$ satisfying $Lf=0$ such that 
   $f(t,y)\to h_\pm(v) \in C^\infty(\bS^2)$ and $f_t,f_{tt}\to 0$ as 
   $t\to\pm\infty$ uniformly for $y$, 
  \item[3.] {\bf [equivalence classes of monopoles]} \\ 
    gauge equivalence classes of 
    monopoles $[(V,A)]$ such that 
    $V(t,y), V_t(t,y)\to 0$ as $t\to\pm\infty$ uniformly for $y$
\end{itemize}
\end{Thm}

\begin{proof}
 By Theorem \ref{Thm:Inverse}, 
 the correspondence {\it 1 $\Leftrightarrow$ 2} is obtained 
 by putting $f:=Rh$ or $h:=h_+$. 
 On the other hand, {\it 2 $\Rightarrow$ 3} is obtained by putting 
 \begin{equation} \label{eq:(V,Theta)<->f}
   V=\del_t f  \quad \text{and} \quad  
   A=- \check *\, \check d\, f. 
 \end{equation} 

 Now we show {\it 3 $\Rightarrow$ 2}. For any $[(V,A)]$ 
 we can take an element $(V,A)$ in this class satisfying the 
 conditions $1^\circ$ and $2^\circ$ in Proposition \ref{Prop:gauge_fixing}. 
 Then by Proposition \ref{Prop:monopole_potential}, we get unique 
 $f\in C^\infty_*(S^3_1)$ satisfying \eqref{eq:(V,Theta)<->f} and $Lf=0$. 
\end{proof}

\begin{Rem}
In the notations in Theorem \ref{Thm:3objects}, 
the evenness $f\in C^\infty_{{\rm even}\, *}(S^3_1)$ and 
the oddness $V\in C^\infty_{\rm odd}(S^3_1)$ automatically hold
by Theorem \ref{Thm:Inverse}.  
\end{Rem}

\paragraph{Admissible monopoles}
To apply Theorem \ref{Thm:3objects} to the study of self-dual metrics, 
we need to assume additional conditions for $(V,A)$, that is, $V$ is positive and 
$V$ is written as in \eqref{eq:Kamada_compactification}. 
Now we introduce the following notion. 

\begin{Def} \label{Def:admissible_monopole}
 Let $(V,A)$ be a monopole on $S^3_1$. Then $(V,A)$ is called {\it admissible} 
 if and only if the following conditions hold: 
 ($1^\circ$) $A(\del_t)=0$, 
 ($2^\circ$) $\check d\, \check *\, A=0$, and 
 ($3^\circ$) $V>0$ and $V$ satisfies the convergence 
  $V(t,y)\to 1$ and $V_t(t,y)\to 0$ as $t\to\pm\infty$ uniformly for $y$. 
\end{Def}

The following Corollary is obviously deduced from Theorem \ref{Thm:3objects} 
and its proof. 

\begin{Cor}
 There is a natural one-to-one correspondence between 
  \begin{itemize}
   \item smooth functions $h\in C^\infty_*(\bS^2)$ satisfying 
      $|\del_tRh(t,y)|<1$, and 
   \item admissible monopoles $(V,A)$, 
  \end{itemize}
 related by $V=1+\del_tRh$ and $A=-\check *\, \check d\, (Rh)$.  
\end{Cor}

For the condition \eqref{eq:Kamada_compactification}, the following hold. 

\begin{Prop} \label{Prop:compactification_condition}
 Let $(V,A)$ be an admissible monopole. 
 Then the condition \eqref{eq:Kamada_compactification} 
 in Proposition \ref{Prop:Kamada_compactification} is always satisfied. 
 Thus any admissible monopole $(V,A)$ defines an self-dual metric 
 $\bar{g}_{V,A}$ on $\bar{\cM}=S^2\times S^2$ 
 with respect to a suitable orientation. 
\end{Prop}
\begin{proof}
 Let $(V,A)$ be an admissible monopole. 
 If we put $\tilde V:=V-1$, then by Theorem \ref{Thm:3objects}
 there exists a generating function $h\in C^\infty_*(\bS^2)$ such that 
 $\tilde V=\del_t Rh=-Q\Delta_{\bS^2}h$. 
 Since $\tilde V$ is odd, it is enough to check the case of $t\to +\infty$. 
 Using the same spherical coordinate $(\c,\p)$ as in the proof of 
 Lemma \ref{Lem:del_Q}, we can write 
 $$  \tilde V (t_0,y) = -\frac{1}{2\pi} \int_{\Omega_{(t_0,y)}} 
          \Delta_{\bS^2}h \, \w_{\bS^2} 
	= -\frac{1}{2\pi} \int_0^\a \left[ \int_0^{2\pi} 
	  \Delta_{\bS^2}h (u(\c,\p)) d\p \right] \sin\c\,d\c, $$
 where $\a$ is defined by $\cos\a=\tanh t_0$. 
 Since the parameter $\c$ is defined by $\cos \c=\tanh t$, 
 $\c$ depends only on $\k:=e^{-2t}$. So we can put 
 $$ {\mathscr F}(\k,y) := \int_0^{2\pi} \Delta_{\bS^2}h (u(\c(\k),\p)) d\p. $$
 Then we obtain 
 $$ \tilde V(t_0,y) = -\frac{1}{\pi} \int_0^{q_0^2} 
  \frac{{\mathscr F}(\k,y)}{(1+\k)^2} \, d\k $$
 where $q_0=e^{-t_0}$. 
 Hence $\tilde V(t,y)$ is a smooth function depending only on $y$ and 
 $q^2=e^{-2t}$, and satisfies $\lim_{q\to+\infty}\tilde V(t,y)=0$. 
 Therefore $V$ is written as in \eqref{eq:Kamada_compactification}. 
\end{proof}

Later (Corollary \ref{Cor:ASD}), 
we will prove the self-duality of the metric $\bar{g}_{V,A}$ on $S^2\times S^2$ 
in a different way from Tod's or Kamada's method. 
(See \cite{bib:Kamada02,bib:Tod93} or the positive definite case 
 \cite{bib:LeBrun91} for their method.)   
By our method, we can determine the `orientation', that is,  
we fix a certain orientation on $S^2\times S^2$ and 
show that $\bar{g}_{V,A}$ is {\it anti-self-dual} with respect to this orientation. 
Moreover, we will see in Corollary \ref{Cor:Zollfreiness} that 
this metric $\bar{g}_{V,A}$ is {\em Zollfrei}.

\paragraph{Example}

Finally in this section, we deal with examples of monopole solutions 
obtained by Tod \cite{bib:Tod93} and Kamada \cite{bib:Kamada05}. 
Let $\{Y^l_m(y)\}_{|m|\le l}$ be the basis of eigenspace of $\Delta_{\bS^2}$ 
with the eigenvalue $-l(l+1)$ 
(i.e. $Y^l_m(y)\in C^\infty(\bS^2)$ can be taken as the 
 {\it spherical harmonics}). 
Introducing variable $z=\tanh t$, let $P_l(z)$ be the Legendre polynomial of 
degree $l$, and put $Z_l(z):=\del_z P_l(z)$. 
In these notations, Tod's monopole solution $(V,A)$ is given by 
$$ V =1+ \sum_{l\ge 1} \sum_{|m|\le l} c_{lm}Z_l(z)Y^l_m(y), \qquad 
 A = - \sum_{l\ge 1} \sum_{|m|\le l} c_{lm} P_l(z) \, 
 \check *\, \check d\, Y^l_m(y), $$
where $\{c_{lm}\}$ is a finite collection of real constants 
with sufficiently small $|c_{lm}|$. 
We remark that the above solution $V$ is first obtained by Tod, and later Kamada 
obtained the above $V$ again with the description of $A$. 
This monopole solution $(V,A)$ is admissible, and the corresponding 
monopole potential $f\in C^\infty_*(S^3_1)$ and the generating function 
$h\in C^\infty_*(\bS^2)$ are given by 
\begin{equation} \label{eq:Tod's_f&h} 
 f= \sum_{l\ge 1} \sum_{|m|\le l} c_{lm} P_l(z) Y^l_m(y), \qquad 
 h= \sum_{l\ge 1} \sum_{|m|\le l}c_{lm} Y^l_m(y). 
\end{equation}

On the other hand, Kamada constructed another type of monopole solutions 
parametrized by the space of probability measures on the hyperboloids 
$H^3_+\sqcup H^3_-$ in the Minkowski space $\R^4_1$. 
However, Theorem \ref{Thm:3objects} insists that this type of solution 
should be gauge equivalent to the above Tod type 
admissible monopole at least asymptotically. 
Actually, Tod's example densely covers all admissible monopoles 
since any generating function $h\in C^\infty_*(\bS^2)$ can be expanded 
as in the form \eqref{eq:Tod's_f&h}.

\section{Local reduction theory}
\label{Sect:Local_theory}

To construct the twistor correspondence for the self-dual metric 
$\bar{g}_{V,A}$ on $S^2\times S^2$ obtained from an 
admissible monopole $(V,A)$, in this section 
we study $S^1$-bundle $\varpi : M^4 \to X^3$ and 
integrable structures on $X$ and $M$.

\paragraph{Einstein-Weyl 3-space} 

Though we only need the integrable property for the 
Einstein manifold $(S^3_1,g_{S^3_1})$, 
we briefly recall the integrability theorem for general 
three-dimensional torsion-free Einstein-Weyl structures 
since there are no difference between the general case and the 
special case of $S^3_1$ so far as studying local theory. 
For the definition of Einstein-Weyl structure, see 
\cite{bib:Hitchin82,bib:LM08,bib:Nakata09}. 
Here we only need the fact that $(S^3_1,[g_{S^3_1}], \nabla^{S^3_1})$ 
is Einstein-Weyl where $\nabla^{S^3_1}$ is the Levi-Civita connection 
of the indefinite metric $g_{S^3_1}$.

Let $X$ be a real 3-manifold, $[g_X]$ be a conformal structure on $X$ 
of signature $(-++)$. We fix a metric $g_X\in [g_X]$ and a frame 
$\{\ul{E}_1,\ul{E}_2,\ul{E}_3\}$ of $TX$ on an open set 
$\ul{U}\subset X$ so that 
\begin{equation} \label{eq:ONF_X}
g_X(\ul E_j,\ul E_k) = 
   \begin{cases} -1 \quad & j=k=1\\ \ 1 & j=k=2 \ \text{or} \ 3 \\ 
     \ 0 & \text{otherwise}.  \end{cases} 
\end{equation}
Let $\nabla^X$ be a torsion-free connection on $TX$, and 
$\ul\w$ be its connection form with respect to the above frame. 
Suppose that $\nabla^X$ is compatible with $[g_X]$, that is, 
$\ul\w$ is written as 
\begin{equation}
 \ul\w = \begin{pmatrix} 
  \p&\ul\w^1_2&\ul\w^1_3 \\ \ul\w^2_1&\p&\ul\w^2_3 \\ 
  \ul\w^3_1&\ul\w^3_2&\p
  \end{pmatrix}, \qquad\qquad 
  \left\{ \begin{aligned} 
   \ul\w^2_1&=\ul\w^1_2, \\ 
   \ul\w^3_1&=\ul\w^1_3, \\ 
   \ul\w^3_2&=-\ul\w^2_3. \end{aligned} \right. 
\end{equation}

A tangent two plane $\bV\subset T_xX \ (x\in X)$ 
is called a {\it null plane} iff $g_X$ 
degenerates on $\bV$, or equivalently, 
iff $\bV$ is tangent to the null cone of $g_X$. 
We put $\bV(\z):=\Span\inn{\ul\Gm_1(\z)}{\ul\Gm_2(\z)}$
for each $\z\in\R\cup \{\infty\}=\RP^1$ where 
\begin{equation} \label{eq:EW_Gm1,2} 
 \ul\Gm_1(\z) :=  -\ul E_1+\ul E_2+\z \ul E_3, \qquad 
 \ul\Gm_2(\z) :=  \z \ul E_1+\z \ul E_2 - \ul E_3. 
\end{equation}
Then $\bV(\z)$ is a null plane, 
and any null plane is written in this form. 

Now let us define the `bundle of null planes' on $X$ by 
$$ \cW_\R :=\left\{ [a]\in\bP(T^*X) 
  \,|\, g_X(a,a)=0 \right\}. $$
Notice that, for each $[a]\in \cW_{\R,x}=\bP(T_x^*X)$, the tangent plane 
$\ker a\subset T_xX$ is a null plane. 
If we  define a 1-form $\ul\Ga(\z)$ by 
$$ \ul\Ga(\z) := (1+\z^2) \ul E^1 + (1-\z^2) \ul E^2 +2\z \ul E^3 $$
using the dual frame $\{\ul E^i\}$ of $\{\ul E_i\}$, 
then we obtain $\bV(\z)=\ker \ul\Ga(\z)$. 
Hence the map 
$\ul U \times \RP^1 \to \cW_\R|_{\ul U}: (x;\z) \mapsto [\ul\Ga(\z)]_x$ 
gives a local trivialization of $\cW_\R$. 
If we introduce coordinates $\theta\in S^1$ by 
$\z=\tan\frac{\c}{2}$ and  $\w = e^{i\c}\in \U(1)$, 
then we obtain the trivializations
\begin{equation} \label{eq:trivialization2_W}
 \ul{U}\times S^1 \overset{\sim}\longrightarrow \cW_\R|_{\ul U} \ : \ 
 (x;\c) \longmapsto 
   [\ul\Ga]_x= [\ul{E}^1 + \cos\c \, \ul{E}^2 + \sin\c\, \ul{E}^3]_x, 
 \hspace{14mm}
\end{equation} 
\begin{equation} \label{eq:trivialization3_W}
 \ul{U}\times \U(1) \overset{\sim}\longrightarrow \cW_\R|_{\ul U} \ : \ 
 (x;\w) \longmapsto 
   [\ul\Ga]_x= [2\w \ul{E}^1 + (1+\w^2) \, \ul{E}^2 +i (1-\w^2) \ul{E}^3]_x. 
\end{equation}
Let us take an open covering $\{\ul U_\a\}$ of $X$ and 
the trivializations of $\cW_\R$ on each $\ul U_\a$ in the form of 
\eqref{eq:trivialization3_W}. 
Then the transition functions are given by the maps 
$F_{\a\b} : \ul U_\a \cap \ul U_\b \to \Aut(\U(1))$ 
where $\Aut(\U(1))$ is the M\"obius transforms on $\U(1)$. 
If $(X,[g_X])$ is {\it space-time orientable}, these transition functions 
can be taken so that $F_{\a\b} : \ul U_\a \cap \ul U_\b \to \Aut(\bD)$ 
where $$\bD:=\{\w\in \C \mid |\w|\le 1 \}$$ 
and $\Aut(\bD)$ is the holomorphic automorphism on $\bD$. 
Hence if $(X,[g_X])$ is space-time orientable, 
we can define the $\bD$-bundle $\cW_+\to X$ associated with 
the $\U(1)$-bundle $\cW_\R\to X$. 
Notice that we obtain a local trivialization 
$\ul U \times \bD \overset\sim \longrightarrow \cW_+|_{\ul U}$ 
by the same equation as \eqref{eq:trivialization3_W} considering $\w\in \bD$. 
We remark that $\cW_+$ is also defined intrinsically 
as the bundle of {\it complex null planes} 
satisfying an orientation compatibility condition (see \cite{bib:Nakata09}). 
We note that the fiber coordinates $\z$ and $\w$ are related by 
$\z=i\frac{1-\w}{1+\w}$, and the disk $\bD=\{|\w|\le 1\}$ corresponds to 
the upper half plane $\{\z\in \C \mid \Imag \z \ge 0 \}$. 

Since the connection $\nabla^X$ is compatible with $g_X$, 
$\cW_\R$ is equipped with a natural connection which we also denote by 
$\nabla^X$. 
Let $\tilde v\in T_{(x;\z)}\cW_\R$ be the horizontal lift 
of a vector $v\in T_xX$ with respect to $\nabla^X$. 
Then by a direct calculation, we obtain the following lifting formula: 
 \begin{equation} \label{eq:EW_lift}
  \tilde{v}= v+ \frac{1}{2} \( \hasira 
  (1+\z^2)\ul\w^2_3 +(1-\z^2)\ul\w^1_3-2\z \ul\w^1_2 \)(v) \, \pd{\z}  
 \end{equation}
Let $\tilde{\ul\Gm}_j$ $(j=1,2)$ be the tautological lift of $\ul\Gm_j$ 
on $\cW_\R$, i.e. 
$ (\tilde{\ul\Gm}_j)_{(x;\z)}=(\ul\Gm_j(\z)_x)\tilde{} \, , $
where $(\cdot)\tilde{}\, $ is the horizontal lift given by \eqref{eq:EW_lift}. 
We define a two-plane distribution on $\cW_\R$ by 
$\ul\cD:=\Span\inn{\tilde{\ul\Gm}_1}{\tilde{\ul\Gm}_2}. $
The integrability of Einstein-Weyl condition is stated as follows. 
\begin{Prop}
 The pair $([g_X],\nabla^X)$ 
 is Einstein-Weyl iff the two-plane distribution $\ul\cD$ is Frobenius integrable. 
\end{Prop}
See \cite{bib:Nakata09} (Proposition 3.9) for the proof.

\paragraph{Indefinite anti-self-dual 4-space} 
Next we summarize the integrable property for 4-dimensional 
anti-self-dual conformal structure of indefinite signature. 
Let $M$ be a real 4-manifold and $[g_M]$ be a conformal structure on $M$ 
of signature $(--++)$. 
We fix $g_M\in[g_M]$ and a frame 
$\{E_0,E_1,E_2,E_3\}$ of $TM$ on an open set $U\subset M$ so that 
\begin{equation} \label{eq:ONF_M}
g_M(E_j, E_k) = 
   \begin{cases} -1 \quad & j=k=0 \ \text{or} \ 1 \\ 
     \ 1 & j=k=2 \ \text{or} \ 3 \\ 
     \ 0 & \text{otherwise}.  \end{cases} 
\end{equation}
The connection form  $\w$ of the Levi-Civita connection $\nabla$ of $g_M$ 
with respect to the above frame is written as 
\begin{equation} \label{eq:SD_conn_form}
 \w = \begin{pmatrix} 
  0&\w^0_1&\w^0_2&\w^0_3 \\  \w^1_0&0&\w^1_2&\w^1_3 \\ 
  \w^2_0&\w^2_1&0&\w^2_3 \\ \w^3_0&\w^3_1&\w^3_2&0
  \end{pmatrix}, \qquad\qquad 
  \left\{ \begin{aligned} 
   \w^1_0&=-\w^0_1, \\ 
   \w^2_0&=\w^0_2, \\ 
   \w^3_0&=\w^0_3, \end{aligned} \right. \qquad 
  \left\{ \begin{aligned} 
   \w^2_1&=\w^1_2, \\ 
   \w^3_1&=\w^1_3, \\ 
   \w^3_2&=-\w^2_3. \end{aligned} \right. 
\end{equation}

We have the eigenspace decomposition 
$\wedge^2 TM=\wedge_+\oplus\wedge_-$ with respect to the Hodge's operator 
on $M$ where $\wedge_\pm$ is the $\pm 1$-eigenspace. 
Using the above frame $\{E_j\}$, we can write as 
$$ \wedge_+ = \Span\left< \vp_1, \vp_2, \vp_3 \right>, \qquad 
 \left\{ \begin{aligned} 
 \sqrt{2}\, \vp_1 &= E_0\wedge E_1 + E_2\wedge E_3, \\ 
 \sqrt{2}\, \vp_2 &= E_0\wedge E_2 + E_1\wedge E_3, \\ 
 \sqrt{2}\, \vp_3 &= E_0\wedge E_3 - E_1\wedge E_2. 
 \end{aligned} \right. $$
Similarly, we have the decomposition 
$\wedge^2 T^*M= \wedge^+\oplus\wedge^-$, 
and we can take a frame $\{\vp^1,\vp^2,\vp^3\}$ of $\wedge^+$ so that 
$\{\vp_j\}$ and $\{\vp^j\}$ are dual each other. 

The Levi-Civita connection $\nabla$ induces a connection on $\wedge_+$ 
which is also denoted by $\nabla$, and its 
connection form is written as  
\begin{equation} \label{eq:yeta<-omega}
  \y = \begin{pmatrix} 0&\y^1_2&\y^1_3 \\ 
  \y^2_1&0&\y^2_3 \\ \y^3_1&\y^3_2&0 \end{pmatrix}, \qquad \qquad
  \left\{ \begin{aligned} 
  \y^2_1&=\y^1_2=\w^1_2-\w^0_3, \\
  \y^3_1&=\y^1_3=\w^1_3+\w^0_2, \\
  -\y^3_2&=\y^2_3=\w^2_3-\w^0_1. \end{aligned} \right. 
\end{equation}

A tangent two plane $\bV\subset T_xM$ is called an {\it $\a$-plane} iff 
$g_M(\bV,\bV)=\{0\}$ (i.e. $\bV$ is contained in the null cone of $g_M$) 
and $\wedge^2 \bV\subset \wedge_+$. 
We put $\bV(\z)=\Span\inn{\Gm_1(\z)}{\Gm_2(\z)}$ for each 
$\z\in \R\cup\{\infty\}=\RP^1$ where  
\begin{equation} \label{eq:SD_Gm1,2}
 \Gm_1(\z) := -\z E_0-E_1+E_2+\z E_3, \qquad
 \Gm_2(\z) := -E_0+\z E_1+\z E_2 - E_3. 
\end{equation}
Then $\bV(\z)$ is an $\a$-plane, 
and each $\a$-plane is written in this form. 

We define the `bundle of $\a$-planes' on $M$ by 
$$ \cZ_\R =\left\{ [\vp]\in\bP(\wedge^+) \,|\, g(\vp,\vp)=0 \right\}. $$
Notice that for each $[\vp]\in \bP(\wedge^+_x)$, the tangent plane 
$\ker\vp:=\{v\in T_xX\mid i(v)\vp=0\}$ is an $\a$-plane. 
If we define $\Ga(\z)\in \wedge^+$ by 
\begin{equation} \label{eq:Ga}
\Ga(\z)= -(1+\z^2)\vp^1 - (1-\z^2) \vp^2 -2\z \vp^3, 
\end{equation}
then we obtain $\bV(\z) =\ker\Ga(\z)$. 
Hence the map $ U\times \RP^1 \to \cZ_\R|_U : (x;\z) \mapsto [\Ga(\z)]_x$ 
gives a local trivialization of $\cZ_\R$. 
Moreover, if $M$ is space-time orientable, we can define the associated 
disk bundle $\cZ_+\to M$ by a similar method as the case of $\cW_+$ 
(see also \cite{bib:LM05}).

The connection $\nabla$ induces a connection on $\cZ_+$ 
which is also denoted by $\nabla$. 
Let $\tilde v\in T_{(x;\z)}\cZ_\R$ be the horizontal lift of a vector 
$v\in T_xM$ with respect to $\nabla$, then 
\begin{equation} \label{eq:SD_lift}
 \tilde{v}= v+ \frac{1}{2} \( \hasira 
  (1+\z^2)\y^2_3 +(1-\z^2)\y^1_3-2\z \y^1_2 \)(v) \, \pd{\z}. 
\end{equation} 
Let $\tilde\Gm_j$ $(j=1,2)$ be the tautological lift of $\Gm_j$ on $\cZ_\R$, i.e. 
$ (\tilde\Gm_1)_{(x;\z)}=(\Gm_1(\z)_x)\tilde{} \, , $
where $(\cdot)\tilde{}\, $ is the horizontal lift given by \eqref{eq:SD_lift}. 
We define a 2-plane distribution on $\cZ_\R$ by 
$ \cD:=\Span\inn{\tilde\Gm_1}{\tilde\Gm_2}. $ 
We can extend $\tilde\Gm_1$ and $\tilde\Gm_2$ to complex vector fields 
on $\cZ_+$ so that they are holomorphic in $\z$. 
We define a complex 3-plane distribution ${\mathcal E}$ on $\cZ_+$ by 
${\mathcal E}:=\Span \langle \tilde\Gm_1,\tilde\Gm_2,\del_{\bar{\z}}\rangle$. 
Then we obtain ${\mathcal E}\cap\overline{\mathcal E}=\{0\}$ on 
$\cZ_+\backslash\cZ_\R$, hence ${\mathcal E}$ defines an almost complex 
structure on $\cZ_+$ so that ${\mathcal E}$ gives the $(0,1)$-vectors. 

\begin{Prop} \label{Prop:integrability_SD}
 The following conditions are equivalent: 
 \begin{itemize}
 \item the conformal structure $[g]$ is anti-self-dual, 
 \item the two-plane distribution $\cD$ on $\cZ_\R$ is Frobenius integrable. 
 \item the almost complex structure on $\cZ_+\backslash\cZ_\R$ 
   defined by $\mathcal E$ is integrable. 
 \end{itemize}
\end{Prop}
See \cite{bib:LM05} (Proposition 3.5 and 7.1) for the proof.

\paragraph{$S^1$-fibration}
Let $(X,g_X)$ be a pseudo-Riemannian 3-manifold of signature $(-++)$ 
and we apply the above argument for $(X,[g_X],\nabla^X)$ where 
$\nabla^X$ is the Levi-Civita connection of $g_X$. 
We put $M:=S^1\times X$ and let $\varpi:M\to X$ be the projection. 
We fix a solution $(V,A)$ of the monopole equation $*dV=dA$ 
on $X$ where $V$ is a positive function 
and $A$ is a one-form on $X$. 
Then $\Theta=ds+A$ defines a connection on the $S^1$-bundle 
$\varpi : M\to X$ 
where $s\in S^1$ is the fiber coordinate. 
We study the following metric on $M$: 
\begin{equation} \label{eq:g_M} 
 g_M := - V^{-2} \Theta\otimes\Theta + g_X. 
\end{equation}
Notice that $g_M$ is conformally equivalent to 
the metric $g_{V,A}=- V^{-1} \Theta\otimes\Theta + V g_X. $

Let us take a local frame $\{\ul E_1,\ul E_2,\ul E_3\}$ 
of $TX$ on an open set $\ul U\subset X$ so that it satisfies the 
orthonormal condition \eqref{eq:ONF_X} for $g_X$. 
We write as $A=A_1 \ul E^1 + A_2 \ul E^2 + A_3 \ul E^3$. 
We define a local frame $\{E_0,E_1,E_2,E_3\}$ of $TM$ on 
$U:=\varpi^{-1}(\ul U)$ by 
\begin{equation} \label{eq:ONF_for_(V,A)}
  E_0=V\pd{s},\quad E_1=\ul{E}_1-A_1\pd{s},\quad
  E_2=\ul{E}_2-A_2\pd{s},\quad E_3=\ul{E}_3-A_3\pd{s}, 
\end{equation}
then $\{E_j\}$ satisfies the orthonormal condition \eqref{eq:ONF_M} for $g_M$. 
Notice that the dual frame $\{E^j\}$ of $\{E_j\}$ is given by 
$$ E^0=V^{-1}\Theta,\quad E^1=\varpi^*\ul{E}^1, \quad E^2=\varpi^*\ul{E}^2, 
  \quad E^3=\varpi^*\ul{E}^3. $$
Now let us use the same notations as above: 
$\ul\w$, $\w$, $\ul\Gm_j$, $\Gm_j$ and so on.

\begin{Lem}
In the above notations, we obtain the following formulas: 
\begin{equation} \label{eq:omega<-(V,A)}
 \begin{aligned}
 \w^0_1&= -\nu_1E^0 + \frac{1}{2}\nu_3E^2 - \frac{1}{2}\nu_2 E^3, \\
 \w^0_2 &= -\nu_2E^0 - \frac{1}{2}\nu_3E^1 - \frac{1}{2}\nu_1 E^3, \\
 \w^0_3 &= -\nu_3E^0 + \frac{1}{2}\nu_2E^1 + \frac{1}{2}\nu_1 E^2, 
 \end{aligned} \qquad\qquad \begin{aligned}
 \w^1_2 &= \varpi^*\ul{\w}^1_2- \frac{1}{2}\nu_3 E^0, \\
 \w^1_3 &= \varpi^*\ul{\w}^1_3+\frac{1}{2}\nu_2 E^0, \\
 \w^2_3 &= \varpi^*\ul{\w}^2_3-\frac{1}{2}\nu_1 E^0, 
 \end{aligned} 
\end{equation}
 where $\nu_j:=V^{-1}E_jV=V^{-1}\ul{E}_jV$ $(j=1,2,3)$. 
\end{Lem}
\begin{proof}
 By the equation $*\,dV=d A$, we obtain  
 $$ \begin{aligned} 
  dE^0 &= d(V^{-1}\Theta) = -V^{-2}dV\wedge\Theta + V^{-1} dA 
      = (V^{-1}\Theta)\wedge (V^{-1}dV) + V^{-1} * dV \\
	&=  E^0\wedge(\nu_1E^1+\nu_2E^2+\nu_3E^3) 
         + (-\nu_1\,E^2\wedge E^3 - \nu_2\,E^1\wedge E^3 
		                  + \nu_3\,E^1\wedge E^2) \\
	&=\sqrt{2}\, (\nu_1\vp^1+\nu_2\vp^2+\nu_3\vp^3). \end{aligned} $$
 Then the required formulas are deduced by a direct calculation 
 so that $\w$ satisfies the torsion-free condition 
 $dE^j+\sum \w^j_k\wedge E^k=0$ and the symmetry \eqref{eq:SD_conn_form}. 
\end{proof}

\begin{Prop} \label{Prop:the_distribution}
In the above notations, we obtain 
\begin{equation} \label{eq:the_distribution}
 \tilde{\Gm}_1  
    = \tilde{\ul{\Gm}}_1 -\( V\z + A(\ul{\Gm}_1) \)\pd{s}, \qquad 
 \tilde{\Gm}_2 
    =\tilde{\ul{\Gm}}_2 - \( V + A(\ul{\Gm}_2)\)\pd{s}. 
\end{equation}
\end{Prop}
\begin{proof}
The proof is given by a direct calculation. Here we sketch the proof of 
the first formula. 
We have $\Gm_1= - \z E_0 -E_1+E_2+\z E_3=
 \ul{\Gm}_1-\z E_0 -A(\ul{\Gm}_1)\del_s$ by definition. 
By the lifting formula \eqref{eq:SD_lift}, we obtain 
$$ \tilde{\Gm}_1 = \Gm_1+ \frac{1}{2} \( \hasira 
  (1+\z^2)\y^2_3 +(1-\z^2)\y^1_3-2\z \y^1_2 \)(\Gm_1) \pd{\z} $$
Evaluating \eqref{eq:yeta<-omega} and \eqref{eq:omega<-(V,A)}, 
and by the lifting formula \eqref{eq:EW_lift}, 
we obtain 
$ \tilde{\Gm}_1 = \tilde{\ul{\Gm}}_1-\z E_0 - A(\ul{\Gm}_1)\del_s 
 = \tilde{\ul{\Gm}}_1 -\(V\z +  A(\ul{\Gm}_1)\)\del_s $
as required. 
\end{proof}

\begin{Rem}
By the result of P.~E.~Jones and K.~P.~Tod \cite{bib:JT}, 
it is natural to expect that, in the above situation, 
the distribution $\cD=\Span\inn{\tilde\Gm_1}{\tilde\Gm_2}$ 
is integrable if and only if 
$\ul\cD=\Span\inn{\tilde{\ul \Gm}_1}{\tilde{\ul \Gm}_2}$ is integrable, 
or equivalently,  $g_M$ is anti-self-dual if and only if 
$([g_X],\nabla^X)$ is Einstein-Weyl. 
To check this claim directly is, however, very hard. 
In the special case of $(S^3_1,g_{S^3_1})$, we prove the integrability of $\cD$ 
by constructing all the integral surfaces of $\cD$
(Proposition \ref{Prop:vanishing}). 
\end{Rem}

\vspace{1ex}
Finally we see that the projection $\varpi:M\to X$ induces a map 
$\Pi:(\cZ_+,\cZ_\R)\to(\cW_+,\cW_\R)$ 
if $X$ is space-time orientable 
(then $M=S^1\times X$ is also space-time orientable). 
For this, notice that $\varpi$ maps each $\a$-plane to a null plane 
since $\varpi_*(\Gm_j)=\ul\Gm_j$ for $j=1,2$. 
Recall that $\cZ_\R$ and $\cW_\R$ are the spaces of $\a$-planes and 
null planes respectively, hence the natural map $\Pi:\cZ_\R\to\cW_\R$ 
is induced. 
By taking local trivializations as above, 
$\Pi$ is locally described as 
$$ \cZ_\R|_U\simeq U\times \RP^1 \longrightarrow 
 \cW_\R|_{\ul U}\simeq \ul U\times \RP^1 : 
  (s,x;\z) \longmapsto (x;\z). $$
Hence the map $\Pi$ naturally extends to a map $\cZ_+ \to \cW_+$. 
By the formula \eqref{eq:the_distribution}, we obtain 
$\Pi_*(\tilde\Gm_i)=\tilde{\ul \Gm}_i$, hence $\Pi_*\cD=\ul \cD$.

\section{Standard model}
\label{Sect:Standard_model}

In this section we study the twistor correspondence for the standard case, 
that is, the case obtained from the trivial monopole $(V,A)=(1,0)$. 

\paragraph{Twistor correspondence for $S^3_1$} 

Recall that we identify the de Sitter space $(S^3_1,g_{S^3_1})$ with 
the space of oriented small circles on $\bS^2$. 
This identification is naturally arisen from 
the LeBrun-Mason correspondence for Einstein-Weyl 3-manifold 
\cite{bib:LM08,bib:Nakata09}. 
Here we describe this correspondence. 

Let us define submanifolds $\Sigma_u\subset S^3_1$ for 
each $u\in \bS^2$ by 
\begin{equation} \label{eq:null_surf}
  \Sigma_u:= \{ (t,y)\in S^3_1 \mid u\in \del \Omega_{(t,y)} \} 
    = \{ (t,y)\in S^3_1 \mid  u\cdot y=\tanh t  \}. 
\end{equation}
Then $\Sigma_u$ gives a {\it null surface}, i.e. $\Sigma_u$ is tangent 
to a null plane at any point on $\Sigma_u$. 
By the correspondence 
$\Sigma_u\leftrightarrow u$, the sphere $\bS^2$ is identified with 
the space of these null surfaces on $S^3_1$. 

Let us introduce the affine coordinates $\l, \y \in \CP^1$ 
related with $y, u\in\bS^2$ by the stereographic projection 
\begin{equation} \label{eq:stereo}
 \l= \frac{y_2+i y_3}{1+y_1}, \qquad 
 \y= \frac{u_2 +i u_3}{1+u_1}. 
\end{equation}
Then the pair $(t,\l)\in \R\times\CP^1$ can be used as the 
coordinate on $S^3_1$. 
We can check by a direct calculation that 
\begin{equation} \label{eq:inequalities}
 u\cdot y > \tanh t \quad \Longleftrightarrow \quad 
  e^t < \left| \frac{\bar{\l}\y+1}{\y-\l} \right|, 
\end{equation}
and that the null surface \eqref{eq:null_surf} is written as 
\begin{equation} \label{eq:null_surf2}
 \Sigma_\y:=\Sigma_u = \left\{ (t,\l)\in \R\times\CP^1 \, \left| \, 
   e^t = \left| \frac{\bar{\l}\y+1}{\y-\l} \right| \right. \right\}. 
\end{equation}

To adapt the formulation in Section \ref{Sect:Local_theory}, we set 
the frame $\{\ul E_j\}$ of $TS^3_1$ on the open set 
$\ul{U}:=\{(t,\l)\in S^3_1 \mid \l\neq\infty\}$ by 
\begin{equation} \label{eq:onf_S^3_1}
 \ul{E}_1=\pd{t}, \quad 
 \ul{E}_2=\frac{1+|\l|^2}{2\cosh t} \(\pd{\l}+\pd{\bar{\l}}\), \quad 
 \ul{E}_3= i\, \frac{1+|\l|^2}{2\cosh t} \(\pd{\l}-\pd{\bar{\l}}\).
\end{equation}
Notice that $\{\ul E_j\}$ satisfies 
the orthonormal condition \eqref{eq:ONF_X} for the metric $g_{S^3_1}$. 
Then the dual frame $\{\ul E^j\}$ is given by 
\begin{equation} 
 \ul E^1 = dt, \quad \ul E^2= \frac{\cosh t}{1+|\l|^2} (d\l+ d\bar{\l}), 
 \quad \ul E^3= -i\, \frac{\cosh t}{1+|\l|^2} (d\l- d\bar{\l}). 
\end{equation}
and the trivialization \eqref{eq:trivialization2_W} is written as 
\begin{equation} \label{eq:(t,l,c)to[a]}
 \ul U\times S^1\overset\sim\longrightarrow \cW_\R|_{\ul U}: 
  (t,\l;\c)\longmapsto [\Ga]= \left[ dt + \frac{\cosh t}{1+|\l|^2} 
   \( e^{-i\c}d\l + e^{i\c}d\bar\l \)\right]. 
\end{equation}
Recall that each point $[\Ga]\in \cW_\R|_x=\bP(T_x^*S^3_1)$ 
corresponds to the null plane $(\ker\Ga)\subset T_xS^3_1$. 
\begin{Prop} \label{Prop:theta<->yeta}
 For each $[\Ga]=(t,\l;\c) \in \ul{U}\times S^1\simeq\cW_\R|_{\ul U}$, 
 the corresponding 
 null plane $\ker\Ga$ is tangent to the null surface $\Sigma_\y$ if and only if 
 \begin{equation} \label{eq:theta<->yeta}
  \y= \frac{-e^{i\c}+\l e^t}{\bar{\l}e^{i\c}+e^t}. 
 \end{equation} 
\end{Prop}
\begin{proof}
 If we put 
 $$ F:= |\y-\l|^2 e^{2t} - |\bar{\l}\y+1|^2, $$ 
 then we can write as $\Sigma_\y=\{(t,\l)\in S^3_1 \mid F=0\}$. 
 Suppose $(t,\l)\in \Sigma_\y$, then the tangent plane $T_{(t,\l)}\Sigma_\y$ 
 is given by $(\ker dF)\subset T_{(t,\l)}S^3_1$, 
 and by a direct calculation we obtain 
 $$ dF= 2 \, |\bar{\l}\y+1|^2 \cdot \left[ dt 
    -\frac{1+|\y|^2}{(\y-\l)(\l\bar\y+1)} d\l 
	-\frac{1+|\y|^2}{(\bar\y-\bar\l)(\bar\l\y+1)} d\bar\l \right]. $$
 Comparing with \eqref{eq:(t,l,c)to[a]}, 
 we see that the coincidence $(\ker\Ga)=(\ker dF)$ occurs if and only if 
 \eqref{eq:theta<->yeta} holds. 
\end{proof}

We put $\bD:=\{\w\in\C\mid |\w|\le 1\}$. 
Since $S^3_1$ is space-time orientable, we can define the 
$\bD$-bundle $\cW_+$ associated with $\cW_\R$, 
and $(t,\l;\w)\in \ul U\times \bD$ gives a local coordinate on $\cW_+|_{\ul U}$. 
We define a smooth map 
$\ul\Gf : \cW_+\rightarrow \CP^1\times\CP^1$ by 
\begin{equation} \label{eq:explicit_ul{Gf}}
  \ul \Gf : (t,\l;\w) \longmapsto (\y_1,\y_2)=
  \(\frac{-\w +\l e^t}{\bar{\l}\w+ e^t}, 
   \frac{\l+e^t\w}{-1+\bar{\l}e^t\w}\) 
\end{equation}
on $\cW_+|_{\ul U}$. 
Then we obtain the double fibration 
\begin{equation} \label{eq:double_fibration_EW}
 \xymatrix{ & (\cW_+,\cW_\R) \ar[dl]_{\ul p} \ar[dr]^{\ul \Gf} \\ 
  S^3_1 && (W,W_\R),}
\end{equation}
where $(W,W_\R)=(\CP^1\times\CP^1,\CP^1)$ and 
$W_\R \hookrightarrow W$ is given by 
$\y \mapsto (\y,\bar{\y}^{-1})$. 
By construction, $\Sigma_\y=\ul p(\ul f^{-1}(\y))$ for each 
$\y\in W_\R\simeq \CP^1$. 
Notice that if we put $\ul D_{(t,\l)} := \ul f (\ul p^{-1}(t,\l))$ then 
$\{\ul D_{(t,\l)}\}_{(t,\l)\in S^3_1}$ 
gives a family of holomorphic disks on $W$ with boundary on 
$W_\R$. 
Further, by the result in \cite{bib:LM08,bib:Nakata09}, the pair 
$([g_{S^3_1}],\nabla^{S^3_1})$ is the unique torsion-free 
Einstein-Weyl structure such that $\{\Sigma_\y\}_{\y\in W_\R}$ gives the 
family of totally geodesic null surfaces on $S^3_1$. 

As easily seen from \eqref{eq:inequalities}, 
the domain $\Omega_{(t,\l)}\subset \bS^2$ coincides with the image of the map 
$$ \bD \longrightarrow \CP^1 : \w \longmapsto 
  \y_1(t,\l;\w)=\frac{-\w+\l e^t}{\bar \l \w + e^t} $$ 
under the identification $\bS^2\overset\sim\to \CP^1$ 
via stereographic projection. 
In particular, the oriented small circle $\del\Omega_{(t,\l)}$ coincides with 
the boundary circle $\del \ul D_{(t,\l)}\subset W_\R$.

\paragraph{Quaternionic description of $S^2\times S^2$} 
Let $\{e_0,e_1,e_2,e_3\}$ be the standard orthonormal basis of 
the Euclidean space $\R^4$, 
and we identify $\R^4$ with the quaternion field $\bH$ by  
$$ a e_0 + b e_1 + c e_2 + d e_3 \ \longleftrightarrow \ 
  a+bi+cj+dk \in \bH.$$

Let $\wedge^2\R^4=\wedge_+\R^4\oplus\wedge_-\R^4$ 
be the eigenspace decomposition for the Hodge's operator on $\R^4$. 
The basis of $\wedge_\pm\R^4$ is given by 
\begin{equation}
  \psi_1^\pm = \frac{1}{\sqrt{2}}\(e_0\wedge e_1 \pm e_2\wedge e_3\), \quad 
  \psi_2^\pm = \frac{1}{\sqrt{2}}\(e_0\wedge e_2 \mp e_1\wedge e_3\), \quad
  \psi_3^\pm = \frac{1}{\sqrt{2}}\(e_0\wedge e_3 \pm e_1\wedge e_2\). 
\end{equation}
Under the identification $\bH\simeq\R^4$, 
we obtain for each $q \in \R^4$ 
\begin{equation} \label{eq:extprod_formula} 
 \begin{aligned}
  \sqrt{2} * ( q \wedge \psi_1^-) &= q i, \\
  \sqrt{2} * ( q \wedge \psi_2^-) &= q j, \\
  \sqrt{2} * ( q \wedge \psi_3^-) &= q k. \\
 \end{aligned}\qquad 
 \begin{aligned}
  \sqrt{2} * ( q \wedge \psi_1^+) &= -i q, \\
  \sqrt{2} * ( q \wedge \psi_2^+) &= -j q, \\
  \sqrt{2} * ( q \wedge \psi_3^+) &= -k q, \\
 \end{aligned} 
\end{equation}
where $*:\wedge^3\R^4\overset\sim\to\R^4$ is the Hodge's operator. 

We define a bilinear form $h$ on $\wedge^2\R^4$ so that it satisfies 
\begin{equation}
 \x_1\wedge\x_2=h(\x_1,\x_2)\, e_0\wedge e_1\wedge e_2\wedge e_3 
\end{equation}
for any $\x_1,\x_2\in \wedge^2\R^4$. 
Then the basis $\{\psi_1^-,\psi_2^-,\psi_3^-,\psi_1^+,\psi_2^+,\psi_3^+\}$ 
gives an orthonormal frame for $h$ of signature $(---+++)$. 
Let us define 
\begin{equation}
 \begin{aligned}
  \cN & := \{ \psi\in\wedge^2\R^4 \mid h(\psi,\psi)=0 
    \ (\text{i.e.}\, \psi\wedge\psi=0)\}, \\
  Q_\R& := \cN/\R_+
   \simeq \left\{ \left. \sum x^l\psi_l^- + \sum y^l\psi_l^+ \, \right| \, 
   x,y\in S^2\right\},
 \end{aligned}
\end{equation}
where the positive real numbers $\R_+$ act on $\cN$ by a scalar multiplication. 
Then $Q_\R$ is diffeomorphic to $S^2\times S^2$ 
and $h$ induces an indefinite conformal structure on $Q_\R$ of signature $(--++)$ 
which is denoted by $[h]$. 
If we define 
\begin{equation}
 \GS_q:=\{ \psi\in Q_\R \mid q \wedge \psi =0 \} 
\end{equation}
for each $q\in S^3\subset\R^4$, then $\GS_q$ 
gives an {\it $\a$-surface} on $(Q_\R,[h])$ 
with respect to the natural orientation on $Q_\R\simeq S^2\times S^2$. 
Since $\GS_q=\GS_{-q}$, the $\a$-surface $\GS_q$ is determined only on 
$[q]\in \RP^3$, so we also write $\GS_{[q]}=\GS_q$.  

By the formula \eqref{eq:extprod_formula}, we can write 
\begin{equation} \label{eq:rep_of_beta_surf}
  \GS_q = \left\{ (x,y)\in S(\Imag\bH)\times S(\Imag\bH) 
    \mid x=\bar{q}yq \right\} 
\end{equation}
under the identification $q\in S^3\simeq \Sp(1)$ and 
$x,y\in S^2\simeq S(\Imag\bH) =\{ \x \in \Imag\bH \mid \x\bar{\x}=1\}$.
If we put $q=a+bi+cj+dk$, then the transform 
$\Imag\bH\to\Imag\bH : y\mapsto \bar{q}yq$ 
is represented by the matrix 
\begin{equation}
 {\mathscr A}(\bar q) := 
 \begin{pmatrix} 
  a^2+b^2-c^2-d^2 & 2(ad+bc) & -2(ac-bd) \\
  -2(ad-bc) & a^2-b^2+c^2-d^2 & 2(ab+cd) \\
  2(ac+bd) & -2(ab-cd) & a^2-b^2-c^2+d^2 
 \end{pmatrix} 
\end{equation}
with respect to the basis $\{i,j,k\}\in\Imag\bH$.
Then we can write as 
$\GS_q=\{ (x,y)\in S^2\times S^2 \mid x={\mathscr A}(\bar q)y \}$. 
We remark that ${\mathscr A}:\Sp(1)\to \SO(3)$ gives a natural double cover. 
By this expression, we see that $\GS_q$ is also an 
$\a$-surface for the standard indefinite metric $g_0$ on $S^2\times S^2$, 
so we obtain $[h]=[g_0]$.

The bundle of $\a$-planes $\hat\cZ_\R\to Q_\R$ is naturally given by 
$$ \hat\cZ_\R = \{ (x,y;[q])\in Q_\R\times\RP^3 \mid (x,y)\in \GS_{[q]} 
 \ \ (\text{i.e.}\ x={\mathscr A}(\bar{q}) y) \}. $$
Since $(Q_\R,[g_0])$ is space-time orientable, 
we can define the disk bundle $\cZ_+$ associated with $\cZ_\R$. 
We will see later (Proposition \ref{Prop:explicit_Gf}) that the projection 
$\hat \Gf : \hat\cZ_\R \to \RP^3$ naturally extends to 
a fiberwise holomorphic map $\hat \Gf : \hat\cZ_+ \to \CP^3$. 
Then we obtain the following double fibration 
(see also \cite{bib:LM05}): 
\begin{equation} \label{eq:double_fibration_SD}
 \xymatrix{ & (\hat\cZ_+,\hat\cZ_\R) \ar[dl]_{\hat p} \ar[dr]^{\hat \Gf} \\ 
  Q_\R && (\CP^3,\RP^3).}
\end{equation}
By construction, we have $\hat p(\hat \Gf^{-1}([q]))=\GS_{[q]}$ 
for each $[q]\in \RP^3$. 
In this way we obtain the LeBrun-Mason twistor space $(\CP^3,\RP^3)$ 
corresponding to the anti-self-dual $4$-manifold $(Q_\R,[g_0])$. 
Here, the two-plane distribution $\cD$ on $\hat\cZ_\R$ is given 
by the tangent distribution of each fiber of $\hat\Gf:\hat\cZ_\R\to \RP^3$.

\paragraph{$S^1$-action} 
Next we study the $S^1$-action on $Q_\R\simeq S^2\times S^2$ defined by 
\eqref{eq:S^1-action}. Recall the notations 
$S_\pm:=\{\pm\e\}\times S^2\subset S^2\times S^2$, 
$\cM:=(S^2\times S^2)\backslash(S_+\sqcup S_-)$, and so on. 
We use the coordinate $(s,t,y)\in S^1\times\R\times \bS^2$ on $\cM$
as in \eqref{eq:x<->(s,t)}. 
Notice that $\cM/S^1=\R\times \bS^2\cong S^3_1$, 
and the quotient map $\varpi: \cM\to S^3_1$ 
is given by $(s,t,y)\mapsto (t,y)$. 
As already mentioned, the standard metric $g_0$ on 
$Q_\R \simeq S^2\times S^2$ is conformally 
equivalent to the metric $g_{1,0}=-ds^2+g_{S^3_1}$ induced from 
the trivial monopole $(V,A)=(1,0)$. 

Now we define the disk bundle $p:(\cZ_+,\cZ_\R)\to \cM$ 
as the restriction of $\hat p : (\hat\cZ_+,\hat\cZ_\R) \to Q_\R$ on $\cM$. 
By the argument in the previous section, 
$\varpi : \cM\to S^3_1$ induces the natural map 
$\Pi:(\cZ_+,\cZ_\R)\to (\cW_+,\cW_\R)$.

Let $\{\ul{E}_1,\ul{E}_2,\ul{E}_3\}$ be the frame of $TS^3_1$ on 
$\ul{U}= \{(t,\l) \mid \l\neq\infty\}$ defined in \eqref{eq:onf_S^3_1}. 
We introduce a frame $\{E_0,E_1,E_2,E_3\}$ 
of $T\cM$ on $U =\varpi^{-1}(\ul{U})=\{(s,t,\l) \mid \l\neq\infty\}$ by 
\eqref{eq:ONF_for_(V,A)}, that is, 
\begin{equation} \label{eq:onf_cM}
 E_0=\pd{s}, \quad E_1=\ul{E}_1, \quad 
 E_2=\ul{E}_2, \quad E_3=\ul{E}_3. 
\end{equation}
If we define $\ul\Gm_j, \Gm_j$ and so on 
similarly as in Section \ref{Sect:Local_theory}, 
we obtain the trivializations $\ul U\times \bD \cong \cW_+|_{\ul U}$ and 
$U\times \bD \cong \cZ_+|_{U}$. 
Hence we can use coordinates $(t,\l;\w)$ on $\cW_+|_{\ul U}$ and 
$(s,t,\l;\w)$ on $\cZ_+|_{U}$. 
In these coordinate, $\Pi:\cZ_+\to\cW_+$ is written as 
$(s,t,\l;\w) \mapsto (t,\l;\w)$. 

The projection $\varpi$ induces a map between the twistor spaces 
in the following way. 
As in \eqref{eq:rep_of_beta_surf}, each $\a$-surface $\GS_q$ 
is defined by the equation $x={\mathscr A}(\bar{q})y$ 
for each $q\in S^3\subset \R^4$. 
In the coordinate $(s,t,y)\in S^1\times\R\times \bS^2 \simeq \cM$, 
this equation is equivalent to the following system: 
\begin{equation} \label{eq:beta_1}
 \frac{e^{is}}{\cosh t}=(1,i,0){\mathscr A}(\bar{q})\cdot y, 
\end{equation}
\begin{equation} \label{eq:beta_2}
 \tanh t=(0,0,1){\mathscr A}(\bar{q})\cdot y. 
\end{equation}
Comparing \eqref{eq:beta_2} with \eqref{eq:null_surf}, 
we see that 
the projection $\varpi$ maps each 
$\a$-surface $\GS_q$ to the null surface 
$\Sigma_u$ where $u=(0,0,1)\cdot{\mathscr A}(\bar{q})$. 
Hence we obtain the natural map 
\begin{equation} \label{eq:pi_intro} 
  \pi: \RP^3\longrightarrow  \bS^2\cong \CP^1 : [q] 
  \longmapsto u =(0,0,1)\cdot{\mathscr A}(\bar{q})
\end{equation}
between the real twistor spaces. 
We will see soon later that 
$\pi$ extends to the map between complex twistor spaces 
and obtain the following commutative diagram: 
\begin{equation} \label{eq:double_fibration_Reduction}
 \xymatrix{ 
   {} & (\cZ_+,\cZ_\R) \ar[dl]_p \ar[dr]^\Gf \ar[d]^\Pi & \\
   \cM \ar[d]_\varpi & (\cW_+,\cW_\R) \ar[dl]_{\ul{p}} \ar[dr]^{\ul{\Gf}} & 
     (Z,Z_\R) \ar[d]^\pi \\ 
   S^3_1 && (W,W_\R). }
\end{equation}
Here $Z_\R=\RP^3\subset\CP^3$ is the standard real submanifold, 
and $Z\subset \CP^3$ is an open set defined later.

\paragraph{Explicit description of the double fibration}
We set an embedding $\bH\cong\R^4\to \C^4$ by 
$$ 	q=a+bi+cj+dk \longmapsto (z_0,z_1,z_2,z_3) $$
\begin{equation} \label{eq:z<->q}
 \begin{aligned} 
 \sqrt{2}\, z_0 &=a -i b + c - i d, \\ 
 \sqrt{2}\, z_1 &=ia - b -ic + d, 
 \end{aligned} \qquad \begin{aligned}
 \sqrt{2}\, z_2 &=-i a - b +i c + d, \\ 
 \sqrt{2}\, z_3 &= a +i b + c +i d. 
 \end{aligned}
\end{equation}
Notice that the image of above embedding is 
$\{ z_3=\bar z_0, z_2= \bar z_1\}$ and 
the image of $\Sp(1)=\{ q\in \bH \mid q\bar q=1 \}$ 
is the set of $(z_i)$ satisfying 
\begin{equation} \label{eq:z_in_Sp(1)}
 z_3=\bar z_0, \quad z_2= \bar z_1, \quad |z_0|^2+|z_1|^2=|z_2|^2+|z_3|^2=1. 
\end{equation}
The above embedding induces the standard embedding 
$\RP^3\to \CP^3$, and we denote its image by 
$Z_\R:=\{[z_i]\in\CP^3 \mid z_3=\bar z_0, z_2= \bar z_1\}.$
Using this notation, the map $\Gf$ is explicitly described in the following way. 

\begin{Prop} \label{Prop:explicit_Gf}
In the above coordinate, the map 
$\Gf : \cZ_\R \to Z_\R : (s,t,\l,\w)\mapsto [z_i]$ is written as 
\begin{equation} \label{eq:explicit_Gf_R}
 \begin{aligned}
  (s,t,\l;\w) \longmapsto
  [z_0:z_1:z_2:z_3]&= \left[ e^{is}\Phi : e^{is}\Phi \, \y : \bar{\y} : 1 \right] 
  \qquad\quad (\w= e^{i\c}), \\
 \text{where} \quad \y &=\frac{-\w+\l e^t}{\bar{\l}\w + e^t}, \quad 
  \Phi= -i \, \frac{\bar\l \w+ e^t}{\l+e^t\w}. 
 \end{aligned}
\end{equation}
 Moreover, the extension $\Gf:\cZ_+\to \CP^3$ is written as 
\begin{equation} \label{eq:explicit_Gf}
  (s,t,\l;\w) \longmapsto
  [z_0:z_1:z_2:z_3]= \left[ e^{is}\Phi : e^{is}\Phi \, \y_1 : \y_2^{-1} : 1 \right] 
  \qquad\quad (|\w|\le 1), 
\end{equation}
where $\y_1$ and $\y_2$ are defined as in \eqref{eq:explicit_ul{Gf}}. 
\end{Prop}

We remark that \eqref{eq:explicit_Gf} can be also written as 
\begin{equation} \label{eq:explicit_Gf2}
  (s,t,\l;\w) \longmapsto \left[
    -ie^{is}(\bar\l\w + e^t) : -ie^{is}(-\w+\l e^t) : 
    -1+\bar\l e^t\w : \l+e^t\w \right].
\end{equation}
We can check by the description \eqref{eq:explicit_Gf} or 
\eqref{eq:explicit_Gf2} that the map 
$\Gf: \cZ_+\to \CP^3$ naturally extends to a smooth map 
$\hat\Gf: \hat\cZ_+\to \CP^3$. 
(The explicit description of $\hat \Gf$ near $S_\pm=Q_\R\backslash \cM$ 
is given in the proof of Proposition \ref{Prop:cM<->Z}. )

\begin{proof}[Proof of \ref{Prop:explicit_Gf}]
Recall that the map $\ul\Gf : \cW_\R \to W_\R : (t,\l;\w)\mapsto \y$ is 
written as 
$$ \y= \y(t,\l;\w)= \frac{-\w+\l e^t}{\bar\l\w+e^t} \qquad (\w=e^{i\c})$$
as in Proposition \ref{Prop:theta<->yeta}. 
Also recall that the map 
$\pi:\RP^3\cong Z_\R\to W_\R\cong S^2: [q]\mapsto u$ 
is given by $u=(0,0,1){\mathscr A}(\bar{q})$. 
Since $\y\in\CP^1$ and $u\in S^2$ are related by the stereographic projection, 
we obtain 
\begin{equation} \label{eq:yeta=z_1/z_0}
 \y=\frac{u_2+i u_3}{1+u_1}
     = \frac{ia- b - ic + d}{a-i b+c-i d} = \frac{z_1}{z_0}, 
\end{equation}
where $(z_i)\in \C^4$ is the image of $q\in \Sp(1)$. 
On the other hand, from the equation \eqref{eq:beta_1}, we obtain 
\begin{equation} \label{eq:alpha_with_z}
\frac{e^{is}}{\cosh t} = i \, \frac{(z_0+\bar\l z_1)(\l z_0- z_1)}{2(1+|\l|^2)}. 
\end{equation}
 
By conditions \eqref{eq:z_in_Sp(1)} and \eqref{eq:yeta=z_1/z_0}, 
there exists $c\in S^1$ satisfying  
\begin{equation} \label{eq:inverse_of_(yeta)}
 (z_0,z_1,z_2,z_3) = \frac{2}{\sqrt{1+|\y|^2}}
   \( e^{\frac{ic}{2}},\y e^{\frac{ic}{2}}, 
       \bar\y e^{-\frac{ic}{2}}, e^{-\frac{ic}{2}}\).  
\end{equation}
Evaluating \eqref{eq:inverse_of_(yeta)} to \eqref{eq:alpha_with_z}, 
we obtain $e^{ic}= e^{is} \Phi$. 
Evaluating this to \eqref{eq:inverse_of_(yeta)} again, 
we obtain the required description \eqref{eq:explicit_Gf_R} of 
$\Gf:\cZ_\R\to Z_\R$. 

The description \eqref{eq:explicit_Gf}
is soon obtained so that the extended map $\Gf:\cZ_+\to Z$ 
is holomorphic in $\w\in \bD$. 
\end{proof}

We need the following Lemma in Section \ref{Sect:Twistor_correspondence}. 
\begin{Lem} \label{Lem:GmPhi}
 Considering $\Phi=\Phi(t,\l;\w)$ as a function on $\cW_+$ or on $\cW_\R$, 
 we obtain $\tilde{\ul \Gm}_1 \Phi = i\z\Phi$ and 
 $\tilde{\ul \Gm}_2 \Phi =i\Phi$. 
\end{Lem}
\begin{proof}
It is enough to check on $\cW_\R$. 
Recall that the distribution $\cD=\Span\inn{\tilde\Gm_1}{\tilde\Gm_2}$
on $\cZ_\R$ is tangent to each fiber of $\Gf:\cZ_\R\to Z_\R$. 
Thus we obtain $\tilde \Gm_j (e^{is}\Phi) =0$ for $j=1,2$ 
by the explicit description \eqref{eq:explicit_Gf_R} of $\Gf$. 
Then by the formula \eqref{eq:the_distribution}, we obtain the 
required equations since we are now studying the case of the 
trivial monopole $(V,A)=(1,0)$. 
\end{proof}

We remark that, since the distribution 
$\ul \cD= \Span\inn{\tilde{\ul\Gm}_1}{\tilde{\ul\Gm}_2}$ 
on $\cW_\R$ is tangent to 
each fiber of $\ul \Gf:\cW_\R\to W_\R$, 
we obtain $\tilde{\ul\Gm}_j\y=\tilde{\ul\Gm}_j\bar\y=0$ 
for $j=1,2$ on $\cW_\R$ where $\y=\y(t,\l;\w)$ is defined above. 
Extending holomorphically, we also obtain 
$\tilde{\ul\Gm}_j\y_k=0$ for $j=1,2$ and $k=1,2$ on $\cW_+$.

\paragraph{The twistor space}
Now let us define an open set $Z\subset \CP^3$ by 
$$ Z := \CP^3 \setminus (L_+ \sqcup L_-) \qquad 
 \text{where} \qquad \left\{ \begin{aligned} 
  L_+ &=\{[z_i]\in\CP^3 \mid z_2=z_3=0\}, \\   
  L_- &=\{[z_i]\in\CP^3 \mid z_0=z_1=0\}. 
  \end{aligned} \right. $$
Further, let us define the holomorphic map 
$\pi :Z\to W=\CP^1\times\CP^1$ by 
$$ \pi :  [z_0:z_1:z_2:z_3]\longmapsto 
 (\y_1,\y_2) = \(\frac{z_1}{z_0},\frac{z_3}{z_2}\). 
$$
Recall that we defined 
$Z_\R := \{[z_i]\in \CP^3 \mid z_3=\bar z_0, z_2=\bar z_1 \}. $ 
Since $\pi(Z_\R)= \{(\y_1,\y_2)\in W \mid \y_2=\bar \y_1^{-1}\} =W_\R$, 
we obtain the map $\pi : (Z,Z_\R) \rightarrow (W,W_\R). $
Notice that this definition of $\pi$ agrees with the above definition of 
$\pi: Z_\R \to W_\R$ in \eqref{eq:pi_intro} or \eqref{eq:yeta=z_1/z_0}. 

The set $Z$ is also obtained in the following way. 
Recall that the $S^1$-action on $\cM$ is written as 
$\a : (s,t,\l)\mapsto (s+\a,t,\l)$. 
Then by \eqref{eq:explicit_Gf_R}, the natural $S^1$-action on $Z_\R$ is induced 
and is written as 
$\a:[z_0:z_1:z_2:z_3]\mapsto [e^{i\a}z_0:e^{i\a}z_1:z_2:z_3].$ 
This $S^1$-action naturally extends to the holomorphic $\C^*$-action on $\CP^3$ 
given by 
\begin{equation} \label{eq:induced_S^1-action}
 \mu\cdot [z_0:z_1:z_2:z_3]\longmapsto
  [\mu z_0:\mu z_1:z_2:z_3] \qquad (\mu\in \C^*). 
\end{equation}
Then $L_+\sqcup L_-$ is just the fixed point set and $Z$ is the free part 
of this action. 
Notice that the map $\pi: Z\to W$ is nothing but the quotient map 
of the above $\C^*$-action. 

By the description \eqref{eq:explicit_Gf2}, we find that the image of 
$\Gf:\cZ_+\to \CP^3$ is contained in $Z$. 
In this way we have obtained the commutative diagram 
\eqref{eq:double_fibration_Reduction}.

\paragraph{Holomorphic disks}
We have already defined the holomorphic disk 
$\ul D_{(t,\l)} := \ul\Gf(\ul p^{-1} (t,\l))$
for each $(t,\l)\in\R\times\CP^1\simeq S^3_1$. 
Similarly, on the diagram \eqref{eq:double_fibration_SD} we put 
$$D_\x:=\hat\Gf (\hat p^{-1}(\x))$$
for each $\x \in Q_\R$. 
Then $\{D_\x\}_{\x\in Q_\R}$ gives a family 
of holomorphic disks on $\CP^3$ with boundaries on $Z_\R$. 
Recall that we defined $Z:=\CP^3\backslash(L_+\sqcup L_-)$. 

\begin{Prop} \label{Prop:cM<->Z}
 The point $\x\in Q_\R$ is contained in $\cM$ 
 if and only if $D_\x \subset Z$. 
 Further, if $\x\in \cM$ then $\pi(D_\x)=\ul D_{\varpi(\x)}$. 
\end{Prop}
\begin{proof}
We change the variable $(s,t,\l)\in S^1\times\R\times\CP^1\cong\cM$ 
to $(\a,\l)\in\C\times\CP^1$ by setting $\a=e^{t+is}$. 
Then $(\a,\l)$ gives a coordinate on an open neighborhood 
$O$ of $S_-$ where $S_\pm=\{\pm\e\}\times S^2$. 
Notice that $S_-=\{(\a,\l)\in O \mid \a=0\}$ and $(\a,\l)=(0,\l)$ 
corresponds to the point $(-\e,\l)$. 

Now recall that $\Gf: \cZ_+\to Z$ is explicitly written as \eqref{eq:explicit_Gf}.
Let us introduce a variable $\w':=e^{is}\Phi(t,\l;\w)$. 
If $e^t < |\l|$, then $\w\mapsto \Phi(t,\l;\w)$ defines an automorphism 
on $\bD$. Hence we can assume $\w'\in \bD$ on $O$ by shrinking $O$ 
if needed. 
Then the triple $(\a,\l;\w')$ gives a 
local coordinate on $\hat \cZ_+|_O$. 
We obtain that the map 
$\hat\Gf : \hat\cZ_+|_O\to \CP^3$ is written as 
$(\a,\l;\w') \longmapsto \left[ \w' : \w' \y_1 : \y_2^{-1} : 1 \right] $ where 
\begin{equation} \label{eq:yeta_from_alpha}
 \y_1=\frac{(1+|\a|^2)\l\w' + i(1+|\l|^2)\a}{(|\a|^2-|\l|^2)\w'}, 
 \qquad \y_2=\frac{|\a|^2-|\l|^2}{-i(1+|\l|^2)\bar\a\w'+(1+|\a|^2)\bar\l}. 
\end{equation} 
Evaluating $\a=0$, we obtain that the disk $D_{(-\e,\l)}$ is given by 
\begin{equation} \label{eq:special_disk-}
  \w' \longmapsto \left[  \w' : -\w'\bar\l^{-1} : -\l^{-1} : 1 \right]
   \qquad (|\w'|\le 1). 
\end{equation}

By a similar argument, the disk $D_{(\e,\l)}$ is given by 
\begin{equation} \label{eq:special_disk+}
  \w' \longmapsto \left[ 1 : \l : \w'\bar{\l} : \w' \right]
  \qquad (|\w|\le 1). 
\end{equation}
Hence each disk $D_{(\pm\e,\l)}$ intersects with $L_+$ or $L_-$, 
so we obtain $D_{(\pm\e,\l)}\not \subset Z$. 

On the other hand, we have $D_\x\subset Z$ for any 
$\x\in \cM$ since the image of $\Gf:\cZ_+\to\CP^3$ is contained in $Z$. 
Hence $\x\in \cM$ if and only if $D_\x\subset Z$. 
The rest statement is obvious by the description \eqref{eq:explicit_Gf}. 
\end{proof}

\paragraph{Compactification of $S^3_1$}
To study the geometry on $S^3_1$, it is convenient to 
consider its compactification. 
Such a picture is actually significant in the study of LeBrun-Mason correspondence 
for Einstein-Weyl 3-manifold (see \cite{bib:LM08,bib:Nakata09}). 

Let $\widehat S^3_1:=Q_\R / S^1$ be the quotient space, 
and $\hat\varpi: Q_\R \to \widehat S^3_1$ be the quotient map. 
Let us write $\hat\varpi(\pm\e,y) = (\pm\infty,y) \in\widehat S^3_1$. 
Then $\widehat S^3_1\simeq [-\infty,+\infty]\times S^2$ 
is considered as the natural compactification of $S^3_1\simeq \R\times S^2$ 
where $[-\infty,+\infty]$ is the natural compactification of $\R$ 
with extra two points $\pm\infty$. 

If we take the limit $t\to\pm\infty$ for the disks $\ul D_{(t,\l)}$ 
on $(W,W_\R)$, then we obtain not a disk but a {\it marked} $\CP^1$. 
Actually by \eqref{eq:explicit_ul{Gf}}, if we put 
$\ul D_{(+\infty,\l)}:=\lim_{t\to+\infty} \ul D_{(t,\l)}$, then 
$\ul D_{(+\infty,\l)}$ is given by 
$$ \{\l\} \times \CP^1 \subset W. $$
In this limit, the boundaries $\del \ul D_{(t,\l)}$ shrink to the point 
$P_{(+\infty,\l)} :=(\l,\bar{\l}^{-1}) \in W_\R$ 
which is considered as 
the marking point of $\ul D_{(+\infty,\l)}$. 
Similarly, $\ul D_{(-\infty,\l)}:=\lim_{t\to-\infty} \ul D_{(t,\l)}$ is given by 
$$ \CP^1 \times \{-\l\}  \subset W $$
equipped with the marking point at 
$P_{(-\infty,\l)}:= (-\bar{\l}^{-1},-\l).$
Notice that by \eqref{eq:special_disk-} or \eqref{eq:special_disk+}, 
we obtain $\pi(D_{(\pm\e,\l)}\cap Z)= P_{(\pm\infty,\l)}$.

Now let us define the maps 
$\chi_\pm : W_\R \to \del\widehat S^3_1$ by 
$\chi_\pm(P_{(\pm\infty,\l)})= (\pm\infty,\l)$. 
Then we can check that 
$\widehat\Sigma_\y= \Sigma_\y \sqcup \{ \chi_+(\y),\chi_-(\y) \}$ 
for each $\y\in W_\R$ 
where $\widehat\Sigma_\y$ is the compactification of 
the null surface $\Sigma_\y$ in $\widehat S^3_1$. 
Similarly, if we put  
$\GC_{(\y_1,\y_2)} := \ul{p}(\ul \Gf^{-1}(\y_1,\y_2))$ 
for $(\y_1,\y_2)\in W\backslash W_\R$, then we obtain 
$\widehat\GC_{(\y_1,\y_2)}=\GC_{(\y_1,\y_2)} \sqcup 
  \{\chi_+(\y_1), \chi_-(\y_2)\}$ 
where $\widehat\GC_{(\y_1,\y_2)}$ is the compactification 
of $\GC_{(\y_1,\y_2)}$. 
We remark that $\GC_{(\y_1,\y_2)}\simeq \R$ is a 
{\it time-like geodesic} on $S^3_1$ (see \cite{bib:Nakata09}).

Finally we remark that, 
in the picture of the correspondence 
$S^3_1\ni (t,y) \leftrightarrow \Omega_{(t,y)}\subset \bS^2$, 
the limit $\lim_{t\to+\infty}\Omega_{(t,y)}$ shrinks to a point $y\in \bS^2$ 
while $\lim_{t\to-\infty}\Omega_{(t,y)}$ wraps the whole $\bS^2$ and 
closes at the point $y\in\bS^2$.

\section{Twistor correspondence}
\label{Sect:Twistor_correspondence}

\paragraph{Main theorem}
In Section \ref{Sect:Standard_model}, we put 
$Z:=\CP^3\backslash(L_+\sqcup L_-)$ and 
$Z_\R:=\{[z_0:z_1:z_2:z_3]\in \CP^3 \mid z_3=\bar z_0, z_2=\bar z_1 \}$, 
and showed the correspondence between 
the map $\pi:(Z,Z_\R)\to (W,W_\R)$ and 
the $S^1$-bundle $\varpi:\cM\to S^3_1$ equipped with the 
standard metrics. 

We now define the deformation of the real submanifold $Z_\R$ in $Z$ by 
\begin{equation} \label{eq:P}
  P_h :=
     \left\{ [z_0:z_1:z_2:z_3]\in Z \ \left| \ z_3=\bar z_0 e^{-h(z_1/z_0)}, \ 
      z_2=\bar z_1 e^{-h(z_1/z_0)} \right. \right\}
\end{equation}
where $h$ is a smooth function on $\CP^1\cong \bS^2$. 
Notice that $P_h=Z_\R$ if $h\equiv 0$, and 
that $P_h$ is invariant under the $\U(1)$-action on $Z$ which is 
defined as the restriction of the $\C^*$-action \eqref{eq:induced_S^1-action}. 
For any real constant $c$, the holomorphic automorphism 
$$ \CP^3 \longrightarrow \CP^3 : [z_0:z_1:z_2:z_3]
  \longmapsto [z_0:z_1:z_2e^c:z_3e^c]$$ 
maps $P_h$ to $P_{h+c}$, so $P_h$ depends on $h$ essentially up to constant. 
So we assume $h\in C^\infty_*(\bS^2)$. 
Then our goal is the following. 

\begin{Thm} \label{Thm:correspondence}
 Let  $(V,A)$ be an admissible monopole, and $h\in C^\infty_*(\bS^2)$ 
 be the corresponding generating function. 
 Then the self-dual metric on $S^2\times S^2$ induced by 
 $(V,A)$ is Zollfrei, and its LeBrun-Mason twistor space is given by 
 $(\CP^3,P_h)$. 
\end{Thm}

\paragraph{Holomorphic disks}

To prove Theorem \ref{Thm:correspondence}, we first construct 
the family of holomorphic disks, and we recover the 
$S^1$-bundle $\varpi : \cM\to S^3_1$. 
Recall that 
for each $(t,\l)\in S^3_1$, the corresponding holomorphic disk on $W$ with 
boundary on $W_\R$ is given by $\ul D_{(t,\l)}=\ul\Gf(\ul p^{-1}(t,\l))$. 

\begin{Prop} \label{Prop:deformed_holomorphic_discs}
There is a unique family of holomorphic disks $\{\GD_{(s,t,\l)}\}$ on $Z$ 
with boundaries on $P_h$ smoothly 
parametrized by $(s,t,\l)\in S^1\times\R\times\CP^1$ 
and satisfying the condition: $\GD_{(s,t,\l)}$ is mapped biholomorphically 
onto $\ul D_{(t,\l)}$ by $\pi:Z\to W$. 
\end{Prop}

If we put $\cM:=S^1\times\R\times\CP^1$, 
we obtain the $S^1$-bundle $\varpi:\cM\to S^3_1: (s,t,\l)\mapsto (t,\l)$ 
by the above Proposition. 

\begin{proof}[Proof of \ref{Prop:deformed_holomorphic_discs}]
First notice that the boundary $\del\ul D_{(t,\l)}$ is given by the image of 
the map  
\begin{equation} 
 \ul \iota : \del \bD \longrightarrow W_\R : 
 \w \longmapsto \y=\y(\w)=\frac{-\w+\l e^t}{\bar\l\w+e^t} 
 \qquad (\w=e^{i\c}\in\del\bD). 
\end{equation}
Then any smooth map $\iota:\del\bD\to P_h$ satisfying 
$\pi\action\iota=\ul \iota$ is written as 
$$ \iota : \w \longmapsto [z_0:z_1:z_2:z_3]=
   \left[ e^{h(\y(\w))}K(\w) : e^{h(\y(\w))}K(\w) \y(\w) : 
     \overline{\y(\w)} : 1 \, \right] $$ 
using a $\U(1)$-valued smooth function $K$ on $\del\bD$. 

Next we deduce the condition for $K$ so that the map $\iota$ 
extends to a holomorphic map on the disk $\bD=\{|\w|\le 1\}$. 
Let us put $H(t,\l;\w) :=h(\y(\w))$ and let 
$$ H(t,\l;\w)= \sum_{k=-\infty}^\infty H_k(t,\l) \w^k $$
be the Fourier expansion. 
We put 
$$ H_+(t,\l;\w) :=\sum_{k>0} H_k(t,\l)\w^k, \qquad 
 H_-(t,\l;\w) :=\sum_{k<0} H_k(t,\l)\w^k. $$
If $\iota$ extends to a holomorphic map on $\bD$, $K$ must be the form 
$$ K(\w) = e^{(H_+(\w) - H_-(\w))} \tilde K(\w) $$ 
where $\tilde K(\w)$ is a holomorphic function on $\bD$ such that 
$\tilde K(e^{i\c})\in \U(1)$. 
Then $\iota$ is written as 
$$ \begin{aligned}
 \iota : \w \longmapsto [z_0:z_1:z_2:z_3] &=
    \left[ e^{2H_+ +H_0}\tilde K(\w) : e^{2H_+ +H_0}\tilde K(\w) \y_1(\w) : 
     (\y_2(\w))^{-1} : 1 \, \right], \\
 \text{where} &\qquad \y_1=\frac{-\w +\l e^t}{\bar{\l}\w+ e^t}, \ \ 
   \y_2=\frac{\l+e^t\w}{-1+\bar{\l}e^t\w}. 
\end{aligned} $$    

If the image of $\iota$ is contained in $Z$, 
then (i) $\tilde K(\w)$ has unique zero on $\bD$ exactly at the pole of 
$\y_1(\w)$, 
and (ii) $\tilde K(\w)$ has unique pole on $\bD$ exactly at the pole of 
$\y_2(\w)$. 
Hence $\tilde K(\w)$ is written as, using a constant $s\in S^1$,
$$ \tilde K(\w) =e^{is}\Phi(\w) \qquad 
  \text{where} \quad \Phi(\w)= -i \, \frac{\bar\l\w+e^t}{\l+e^t\w}. $$
Thus $\iota$  is written as 
\begin{equation} \label{eq:deformed_disk}
 \iota : \w \longmapsto [z_0:z_1:z_2:z_3]=
   \left[ e^{2H_+ +H_0+is}\Phi : e^{2H_+ +H_0+is}\Phi \y_1 : 
     \y_2^{-1} : 1 \, \right]. 
\end{equation}
Let us define $\GD_{(s,t,\l)}$ to be the holomorphic disk obtained by 
\eqref{eq:deformed_disk}. 
Then the statement follows 
since $H_+,H_0$ and $\Phi$ depends smoothly on $(t,\l)$, 
and are independent of $s$. 
\end{proof}

Recall that 
the boundary $\del {\ul D}_{(t,\l)} \subset W_\R \simeq\CP^1$  
corresponds to the oriented small circle $\del\Omega_{(t,y)}$. 
Hence, in the above proof, the Fourier coefficient $H_0(t,\l)$ is written as 
\begin{equation}
 H_0(t,\l)= \frac{1}{2\pi} \int_0^{2\pi} H(t,\l;e^{i\c}) d\c 
   = \frac{1}{2\pi} \int_0^{2\pi} h(\y(e^{i\c})) d\c 
   = Rh(t,\l)
\end{equation}
using the transform $R$ defined in \eqref{eq:basic_transform}. 
Here we abused the notations as $Rh(t,y)=Rh(t,\l)$.

\paragraph{Non-admissible deformations}

Let $\{\GD_\x\}_{\x\in\cM}$ be the family of holomorphic disks 
obtained in Proposition \ref{Prop:deformed_holomorphic_discs}. 
Let us denote the interior of the disk $\GD_\x$ by $\GD_\x^\circ$. 
We will see later (Proposition \ref{Prop:properties_of_disks}) that 
the family $\{\GD_\x^\circ\}_{\x\in\cM}$ 
foliates $\pi^{-1}(W\backslash W_\R)\subset Z$ 
if the corresponding monopole $(V,A)$ 
is admissible, that is, if the generating function $h\in C^\infty_*(\bS^2)$ 
satisfies $|\del_t Rh(t,\l)|<1$. 

On the other hand, in the non-admissible case, 
we obtain the following. 

\begin{Prop} \label{Prop:non-admissible_case}
 If $|\del_t Rh(t,\l)|>1$ for some $(t,\l)\in S^3_1$, then the family  
 $\{\GD_\x^\circ\}_{\x\in\cM}$ does not give a foliation. 
\end{Prop}
\begin{proof}
 Suppose that there exists a point $(t_0,\l_0)\in S^3_1$ such that 
 $|\del_t Rh(t_0,\l_0)|>1$. 
 Since $Rh(t,\l)$ is an even function on $S^3_1$, 
 we can assume $\del_tRh(t_0,\l_0)+1<0$ by changing $(t_0,\l_0)$ 
 with $(-t_0,-\l_0)$ if needed.  
 
 Now if we evaluate $\w=0$ to the description \eqref{eq:deformed_disk} 
 of the disk $\GD_{(s,t,\l)}$, we find that the disk $\GD_{(s,t,\l)}$ contains 
 the point 
 $ \left[ -ie^{H_0+t+is}: -ie^{H_0+t+is}\l : -1: \l \right]\in Z\backslash P_h$. 
 We claim that the map $ S^1\times \R \to Z\backslash P_h $ given by  
 \begin{equation} \label{eq:non_injective}
 (s,t) \longmapsto  
 \left[ -ie^{H_0(t,\l_0)+t+is}: -ie^{H_0(t,\l_0)+t+is}\l_0 : -1: \l_0 \right] 
 \end{equation} 
 is not injective. 
 If this map is injective, then the function $H_0(t,\l_0)+t$ must 
 be monotonic in $t\in \R$. We have, however, 
 $$\del_tRh(t_0,\l_0)+1<0 \qquad \text{and} \qquad 
   \lim_{t\to\infty} (\del_tRh(t,\l_0)+1)=1>0, $$ 
 hence the function $H_0(t,\l_0)+t=Rh(t,\l_0)+t$ is not monotonic. 
 So the map \eqref{eq:non_injective} is not injective. 
 This means that some members in $\{\GD_{(s,t,\l_0)}\}_{(s,t)\in S^1\times\R}$ 
 intersect with each other at their interior points, 
 hence $\{\GD_\x^\circ\}_{\x\in\cM}$ does not give a foliation. 
\end{proof}

\paragraph{Double fibration}
Next we construct the double fibration. 
Let $(V,A)$ be the monopole corresponding to $h\in C^\infty_*(\bS^2)$, 
and suppose that $(V,A)$ is admissible. 
By Proposition \ref{Prop:compactification_condition}, 
we obtain an indefinite metric $\bar g_{V,A}$ on $\bar{\cM}=S^2\times S^2$. 
Here we show that this metric is anti-self-dual with respect to 
the natural orientation on $\cM=S^1\times\R\times \CP^1$. 

Let $(\hat\cZ_+,\hat\cZ_\R)$ be the disk bundle on 
$S^2\times S^2$ induced from $\bar g_{V,A}$ by the method explained in 
Section \ref{Sect:Local_theory}. 
Recall that $\hat\cZ_\R$ is equipped with the two-plane distribution 
$\cD$ which is locally written as 
$\cD=\Span\inn{\tilde\Gm_1}{\tilde\Gm_2}$, 
and that $\bar{g}_{V,A}$ is anti-self-dual if and only if $\cD$ is integrable. 

Let $(\cZ_+,\cZ_\R) :=(\hat\cZ_+|_\cM,\hat\cZ_\R|_\cM)$ be the 
restriction on $\cM$. 
We take a local trivialization of $(\cZ_+,\cZ_\R)$ 
on the open set $\ul U:=\{(t,\l)\in S^3_1 \mid \l\neq\infty \}$  
in the following way. 
We fix a frame $\{\ul E_j\}_{j=1,2,3}$ of $TS^3_1$ on $\ul U$
in the same way as \eqref{eq:onf_S^3_1}. 
We define the frame $\{E_j\}_{j=0,1,2,3}$ of $T\cM$ 
on the open set $U:=\varpi^{-1}(\ul U)=\{(s,t,\l)\in \cM \mid \l\neq\infty\}$ 
of $\cM$ by \eqref{eq:ONF_for_(V,A)} 
so that we can apply the argument in Section \ref{Sect:Local_theory}. 
Then we obtain the trivialization $ U\times \bD \cong \cZ_+|_U$, 
and we can use $(s,t,\l;\w) \in U \times \bD$ 
as a local coordinate on $\cZ_+|_U$.

Now let $\{\GD_{(s,t,\l)}\}$ be the holomorphic disks obtained in 
Proposition \ref{Prop:deformed_holomorphic_discs}. 
Noticing the explicit description \eqref{eq:deformed_disk} 
of the disk $\GD_{(s,t,\l)}$, 
we define the map 
$(\cZ_+|_U,\cZ_\R|_U)\to (Z,P_h)$ by 
\begin{equation} \label{eq:explicit_deformed_Gf}
  (s,t,\l;\w) \longmapsto [z_0:z_1:z_2:z_3]=
   \left[ e^{2H_+ +H_0+is}\Phi : e^{2H_+ +H_0+is}\Phi \, \y_1 : \y_2^{-1} 
   : 1 \, \right]. 
\end{equation}
It is checked that this map uniquely extends to a smooth map 
$\Gf: (\cZ_+,\cZ_\R)\to (Z,P_h)$.  
In this way, we obtain the similar diagram as 
\eqref{eq:double_fibration_Reduction}. 
By construction, this diagram commutes. 

\begin{Prop} \label{Prop:vanishing}
 In the above notations, each fiber of the map $\Gf: \cZ_\R\to P_h$ is tangent 
 to the distribution $\cD|_{\cZ_\R}$. 
\end{Prop}
\begin{proof}

By the explicit description \eqref{eq:explicit_deformed_Gf} of the map 
$\Gf:\cZ_\R\to P_h$, it is enough to check that 
the following formulas hold for $j=1,2$:  
\begin{equation} \label{eq:wanted} 
 \tilde{\Gm}_j \( e^{2H_++H_0+is}\Phi \) =0, 
\end{equation}
\begin{equation} \label{eq:Gm_j_yeta_k} 
  \tilde{\Gm}_j \y_1= \tilde{\Gm}_j \y_2= 0, 
\end{equation}
on $\cZ_\R|_U$.
The equation \eqref{eq:Gm_j_yeta_k} is, however, obvious 
since the vectors $\Pi_*(\Gm_j)=\tilde{\ul\Gm}_j$ $(j=1,2)$ and 
the functions $\y_k$ $(k=1,2)$ are not deformed from the standard case, so 
$\tilde\Gm_j \y_k= \tilde{\ul \Gm}_j \y_k = 0$ for each $j,k$. 
On the other hand, by Proposition \ref{Prop:the_distribution}, 
the equations \eqref{eq:wanted} 
is equivalent to the following equations: 
$$ \left\{ \begin{aligned}
  \tilde{\ul \Gm}_1 (e^{2H_+ + H_0}\Phi) & = 
   i (V\z + A(\ul\Gm_1)) \cdot e^{2H_+ + H_0}\Phi, \\ 
  \tilde{\ul \Gm}_2 (e^{2H_+ + H_0}\Phi) & = 
   i (V + A(\ul\Gm_2)) \cdot e^{2H_+ + H_0}\Phi, 
 \end{aligned} \right. $$
where $\z=i\, \frac{1-\w}{1+\w}$. 
If we apply Lemma \ref{Lem:GmPhi}, the wanted equation 
\eqref{eq:wanted} is equivalent to 
\begin{equation} \label{eq:wanted2}
 \left\{ \begin{aligned}
  \tilde{\ul \Gm}_1 (2H_+ + H_0) &= i ((V-1)\z + A(\ul\Gm_1)), \\
  \tilde{\ul \Gm}_2 (2H_+ + H_0) &= i ((V-1) + A(\ul\Gm_2)). 
 \end{aligned} \right. 
\end{equation}

Now notice that for $|\w|=1$ we have 
\begin{equation} \label{eq:ul{Gm}H=0} 
 \ul\Gm_j H(t,\l;\w)= \ul \Gm_j h(\y(\w)) 
 = \PD{h}{\y}\cdot\ul \Gm_j \y_1 + \PD{h}{\bar{\y}}\cdot \ul \Gm_j \y_2 = 0. 
\end{equation}
If we use the formula \eqref{eq:Gm_in_omega} in Appendix \ref{App:Formulas}, 
we obtain the following equations 
\begin{equation} 
 l_j H_k + \bar l_j H_{k-1} + k \d_j H_k - (k-1)\bar \d_j H_{k-1} = 0 
 \qquad\quad (k\in \Z) 
\end{equation}
where $H(t,\l;\w)=\sum_k H_k(t,\l)\w^k$. 
Thus we obtain for $j=1$ 
\begin{equation} \label{eq:Gm1(2H_++H_0)}
 \begin{aligned}
 (1+\w)\tilde{\ul \Gm}_1 (2H_+ +H_0) 
  &= 2 ( l_1 H_1+ \d_1 H_1) \w + (l_1H_0 + \w\, \bar l_1 H_0) 
   = l_1H_0 - \w\, \bar l_1 H_0.  \\
  &= -(1-\w) \ul E_1H_0 + (1-\w) \ul E_2H_0 + i(1+\w) \ul E_3H_0. 
 \end{aligned}
\end{equation}

On the other hand, we have $H_0(t,\l)=Rh(t,\l)$ and by the hypothesis 
$$ \begin{aligned}
 V &= 1+\del_t Rh=1+\ul{E}_1H_0, \\ 
 A &= - \check * \, \check d \, Rh 
   = (\ul E_3 H_0) \ul E^2 - (\ul E_2 H_0) \ul E^3. 
   \end{aligned} $$
Hence 
\begin{equation} \label{eq:wanted_RHS}
 \begin{aligned}
   -(1-\w)(V-1) + i A((1+\w)\ul\Gm_1) 
  &= -(1-\w) \ul E_1H_0 + i A(l_1+\w \bar l_1) \\
  &= -(1-\w) \ul E_1H_0 + (1-\w) \ul E_2H_0 + i (1+\w) \ul E_3H_0. 
 \end{aligned}
\end{equation}
By \eqref{eq:Gm1(2H_++H_0)} and \eqref{eq:wanted_RHS}, we obtain 
$$  (1+\w)\tilde{\ul \Gm}_1 (2H_+ +H_0) 
  = -(1-\w)(V-1) + i A((1+\w)\ul\Gm_1) $$
which is equivalent to the first equation of \eqref{eq:wanted2}. 
The second equation of \eqref{eq:wanted2} is proved in a similar way. 
\end{proof}

\begin{Cor} \label{Cor:ASD}
 Let $(V,A)$ be an admissible monopole. 
 Then the metric $\bar g_{V,A}$ on $S^2\times S^2$ induced from $(V,A)$  
 is anti-self-dual with respect to the natural orientation on 
 $\cM=S^1\times \R \times \CP^1$. 
\end{Cor}
\begin{proof}
 Notice that the map $\Gf:\cZ_\R\to P_h$ is surjective by construction. 
 Hence each fiber of $\Gf:\cZ_\R\to P_h$ is two-dimensional, and is an 
 integral surface of $\cD$ by Proposition \ref{Prop:vanishing}. 
 Thus $\cD$ is Frobenius integrable.  
 Hence $\bar g_{V,A}$ is anti-self-dual on $\cM$ 
 by Proposition \ref{Prop:integrability_SD}. 
 Since $\cM$ is dense in $S^2\times S^2$, 
 $g$ is anti-self-dual on the whole of $S^2\times S^2$. 
\end{proof}

By Proposition \ref{Prop:integrability_SD}, 
the complex 3-plane distribution ${\mathcal E}=\Span \langle 
 \tilde{\Gm}_1, \tilde{\Gm}_2, \del_{\bar{\w}}\rangle$ defines  
the complex structure on $\cZ_+\backslash\cZ_\R$. 
Since $e^{2H_++H_0+is}\Phi$ is holomorphic in $\w\in\bD$, 
the equations \eqref{eq:wanted} and \eqref{eq:Gm_j_yeta_k} 
hold on $\cZ_+$. 
Hence the map $\Gf:\cZ_+\to Z$ is holomorphic on 
$\cZ_+\backslash\cZ_\R$. 
In this way, we have obtained the following result. 

\begin{Prop} \label{Prop:deformed_Gf}
  In the above notations, 
 $\Gf: (\cZ_+\backslash\cZ_\R) \to (Z\backslash P_h)$ is holomorphic. 
\end{Prop}

\paragraph{Compactification}

Recall that the compactification $\cM\hookrightarrow S^2\times S^2$ 
is given by $(s,t,\l) \mapsto (x,y) $
where $y \overset\sim\mapsto \l$ is the stereographic projection and 
$(s,t)\mapsto x$ is given by \eqref{eq:x<->(s,t)}. 
We have $(S^2\times S^2)\backslash \cM= S_+\sqcup S_-$ 
where $S_\pm=\{\pm\e\}\times S^2$. 
Similar as the proof of Proposition \ref{Prop:cM<->Z}, 
let us introduce the variables $\a=e^{t+is}$ and 
$\w'=e^{is}\Phi(t,\l,\w)$, then 
$(\a,\l)$ gives a coordinate on the small open neighborhood 
$O\subset S^2\times S^2$ of $S_-$ and 
$(\a,\l;\w')$ gives a local coordinate on $\hat\cZ_+|_O$. 
The map $\Gf:\cZ_+\to Z$ defined in \eqref{eq:explicit_deformed_Gf} 
is written as 
\begin{equation} \label{eq:deformed_Gf_from_alpha}
 \Gf : (\a,\l;\w')\longmapsto 
 \left[ e^{2H_++H_0}\w' : e^{2H_++H_0}\w'\y_1 : \y_2^{-1} :1 \right] 
\end{equation}
where $\y_1$ and $\y_2$ are given by \eqref{eq:yeta_from_alpha}. 
Since the function $H(t,\l;e^{i\c})=h(\y_1(t,\l;e^{i\c}))$ extends to a smooth 
function on $(\a,\l;\w')\in \hat\cZ_\R|_O$, 
its Fourier coefficient $H_k$ extends to  
a smooth function on $(\a,\l)\in O$ for each $k\in \Z$. 
Hence \eqref{eq:deformed_Gf_from_alpha} extends to the smooth map 
$\hat\Gf:\hat\cZ_+|_O\to \CP^3$. 
By a similar argument for $S_+$, we obtain the smooth map 
$\hat\Gf:\hat\cZ_+\to \CP^3$ as an extension of $\Gf$. 
So we get the double fibration 
\begin{equation} \label{eq:deformed_double_fibration_SD}
 \xymatrix{ & (\hat\cZ_+,\hat\cZ_\R) \ar[dl]_{\hat p} \ar[dr]^{\hat\Gf} \\ 
  S^2\times S^2 && (\CP^3,P_h).}
\end{equation}

Let us define the holomorphic disks $\{\GD_\x\}_{\x\in S^2\times S^2}$ by 
$\GD_\x:= \hat{\Gf}(\hat p^{-1}(\x))$. 
Of course, this notation agrees with the previous notation of $\GD_\x$ 
for $\x\in\cM$. 
Since we have 
$$ H_+(\a,\l;\w')|_{\a=0}=\lim_{t\to -\infty} H_+(t,\l,\w)=0, \quad 
  H_0(\a,\l;\w')|_{\a=0}=\lim_{t\to -\infty} H_0(t,\l,\w)=h(-\bar{\l}^{-1}), $$ 
the disk $\GD_{(-\e,y)}$ is given by the map 
$$ \w' \longmapsto   \left[ e^{h(-\bar\l^{-1})}\w' : 
   -e^{h(-\bar\l^{-1})}\w'\bar\l^{-1} : -\l : 1 \,  \right]. $$ 
Similarly, we can check that $\GD_{(\e,y)}$ is given by the map 
$$ \w' \longmapsto 
  \left[ e^{h(\l)}:e^{h(\l)}\l : \w'\bar\l : \w' \right]. $$

The family of holomorphic disks $\{\GD_\x\}$ on $(\CP^3,P_h)$ 
has the following properties. 

\begin{Prop} \label{Prop:properties_of_disks}
 The family $\{\GD_\x\}_{\x\in S^2\times S^2}$ satisfies the following 
 conditions: 
 \begin{enumerate}
  \item $\x\in\cM$ if and only if $\GD_\x\subset Z$, 
  \item for each disk $\GD_\x$, the class $[\GD_\x]\in H_2(\CP^3,P;\Z)$ 
             gives a generator, and 
  \item $\{\GD_\x^\circ\}_{\x\in S^2\times S^2}$ 
  foliates $\CP^3\backslash P_h$ where $\GD_\x^\circ$ is the interior of 
  $\GD_\x$, 
 \end{enumerate}
\end{Prop}
\begin{proof}
 The statement {\it 1} is easily deduced by the above descriptions of $\GD_{\x}$. 
 To check {\it 2}, 
 it is enough to check the case $\x=(\pm\e,y)\in S_\pm$ since 
 all the disks of $\{\GD_\x\}$ are homotopic in $(\CP^3,P_h)$ each other, 
 and this is obvious by the above descriptions of $\GD_{(\pm\e,y)}$. 

 To prove {\it 3}, we show 
 \begin{enumerate}
 \item[$({\rm 1}^\circ)$] the family $\{\GD_\x^\circ\}_{\x\in\cM}$ foliates 
    $Z^\circ := \pi^{-1}(W\backslash W_\R)$, and 
 \item[$({\rm 2}^\circ)$] the family 
    $\{\GD_\x^\circ\}_{\x\in (S_+\sqcup S_-)}$ foliates 
    $(\CP^3\backslash Z^\circ)\backslash P_h$. 
 \end{enumerate}
 Here ${\rm 2}^\circ$ is obviously deduced by the descriptions of 
 $\GD_{(\pm\e,y)}$ since 
 $$ (\CP^3\backslash Z^\circ)\backslash P_h 
  = \left\{ \left. \left[ce^{h(\y)}: ce^{h(\y)}\y : \bar\y : 1 \right] \in\CP^3 
     \ \right| \ \y, c \in\CP^1, |c|\neq 1 \right\}.$$ 
 To check ${\rm 1}^\circ$, it is enough to show that 
 $\Gf: (\cZ_+\backslash\cZ_\R) \to Z^\circ$ is bijective. 
 For this, we only need to show that the restriction 
 $\Gf^{-1}(\pi^{-1}(\y_1,\y_2)) \overset\Gf \to 
   \pi^{-1}(\y_1,\y_2)$ is bijective 
 for each point $(\y_1,\y_2)\in W\setminus W_\R$. 
 We put ${\mathcal C}:= \Gf^{-1}(\pi^{-1}(\y_1,\y_2)) \subset \cZ_+$ and 
 $\ul{\mathcal C}:= \ul\Gf^{-1}(\y_1,\y_2) \subset \cW_+$. 
 The set $\ul {\mathcal C}$ is diffeomorphic to $\R$. 
 (In fact this is a canonical lift of the time-like geodesic $\GC_{(\y_1,\y_2)}$. 
 See \cite{bib:Nakata09}.)  
 The map $\Pi|_{\mathcal C}:{\mathcal C}\to \ul{\mathcal C}$ is an 
 $S^1$-bundle with fiber coordinate $s\in S^1$. 
 Notice that 
 ${\mathcal C}$ is a complex submanifold 
  of $\cZ_+\backslash\cZ_\R$ since $\Gf$ and $\pi$ are holomorphic. 

 Now we suppose $\y_1\neq\infty$ and $\y_2\neq 0$.  
 First we check that $\Gf|_{\mathcal C}: {\mathcal C} \to \pi^{-1}(\y_1,\y_2)$ 
 is injective. 
 In this case we can take a 
 coordinate $\C^*\cong \pi^{-1}(\y_1,\y_2)$ by 
 $\mu  \mapsto [\mu: \mu \y_1: \y_2^{-1} : 1] $. 
 Then  $\Gf|_{\mathcal C}$ is written as 
 $\cF : (s,t,\l;\w)\mapsto e^{2H_+ +H_0+is}\Phi$ 
 by \eqref{eq:explicit_deformed_Gf}. 
 Notice that $\del_s\cF=i\cF$. If we denote the complex structure on 
 ${\mathcal C}$ by $J$, then $(\del_s+iJ(\del_s))\cF=0$ since 
 $(\del_s+iJ(\del_s))$ is a (0,1)-vector field. Hence $J(\del_s)\cF=-\cF$. 
 Here by $\y_1\neq\infty$ and $\y_2\neq 0$, we obtain $\cF\neq 0,\infty$. 
 So $J(\del_s)\cF=-\cF$ means that any two fiber-circles of 
 ${\mathcal C}\to \ul{\mathcal C}$ are mapped by $\Gf$ to different  
 $\U(1)$-orbits in $\pi^{-1}(\y_1,\y_2)\cong \C^*$. 
 Hence $\cF$ is injective. 

 Next we check the surjectivity of $\cF$. 
 For this, it is enough to show that $\lim_{t\to-\infty} \cF= 0 $ 
 and $\lim_{t\to +\infty}\cF =\infty$. 
 As explained in the last part of Section \ref{Sect:Standard_model}, 
 $\GC_{(\y_1,\y_2)}=\ul p(\ul \Gf^{-1}(\y_1,\y_2))$ 
 is a time-like geodesic connecting 
 $\chi_+(\y_1)$ and $\chi_-(\y_2)$ in the compactification 
 $\widehat S^3_1$ of $S^3_1$. 
 Recall that $\varpi:\cM\to S^3_1$ is naturally compactified to the quotient 
 map $\hat\varpi:S^2\times S^2\to \widehat S^3_1$. 
 So the set $\varpi^{-1}(\GC_{(\y_1,\y_2)})\simeq \R\times S^1$ is 
 compactified to $\hat\varpi^{-1}(\widehat\GC_{(\y_1,\y_2)}) \simeq S^2$ 
 in $S^2\times S^2$ 
 with extra two points $(+\e,\y_1)$ and $(-\e,\y_2)$. 
 Notice that ${\mathcal C}$ is homeomorphically mapped onto 
 $\varpi^{-1}(\GC_{(\y_1,\y_2)})$ by $p$. 
 For any path $\g(\t)$ in ${\mathcal C}$ with parameter $\t\in [0,\infty)$ 
 such that $\lim_{\t\to\infty} p(\g(\t))=(+\e,\y_1)$, 
 we obtain $\lim_{\t\to\infty}\Gf(\g(\t))$ is, if exist, contained in the 
 disk $\GD_{(+\e,\y_1)}$. 
 On the other hand, $\lim_{\t\to\infty}\Gf(\g(\t))$ is, if exist, contained in 
 the closure $\overline{\pi^{-1}(\y_1,\y_2)}$. 
 Since the intersection $\overline{\pi^{-1}(\y_1,\y_2)}\cap\GD_{(+\e,\y_1)}$ 
 is one point $z=[1:\y_1:0:0]$, the limit $\lim_{\t\to\infty}\Gf(\g(\t))$ 
 actually exist independently with the path $\g(\t)$  
 and the limit is above $z$. 
 Hence we obtain $\lim_{t\to\infty}\cF=\infty$. 
 Similarly we can check $\lim_{t\to -\infty}\cF=0$. 
 
 In this way, we have proved that $\Gf|_{\mathcal C}$ is bijective 
 if $\y_1\neq\infty$ and $\y_2\neq 0$. 
 in the case $\y_1=\infty$ or $\y_2=0$, 
 we can check the bijectivity of $\Gf|_{\mathcal C}$ similarly 
 by taking a suitable coordinate $\C^*\cong \pi^{-1}(\y_1,\y_2)$. 
 Thus the statement {\it 3} is proved. 
\end{proof}

\paragraph{Zollfrei condition}
On the double fibration \eqref{eq:deformed_double_fibration_SD}, 
we set 
$$ \GS_z:=\hat p (\hat \Gf^{-1}(z)) 
 = \{ \x\in S^2\times S^2 \mid z\in \GD_\x \} $$
for each point $z\in P_h$. 

\begin{Prop} \label{Prop:closedness_of_a_surf}
 The set $\GS_z$ is a smoothly embedded $S^2\subset S^2\times S^2$  
 for each point $z\in P_h$.
\end{Prop}
\begin{proof}
 We notice to the null surface $\Sigma_{\pi(z)}:=\ul p (\ul \Gf^{-1}(\pi(z)))$ 
 on $S^3_1$. 
 By definition, for each $(t,\l)\in \Sigma_{\pi(z)}$ we have 
 $\pi(z)\in \del \ul D_{(t,\l)}$. 
 Recall that $\{\GD_{(s,t,\l)}\}_{s\in S^1}$ gives the family 
 of disks satisfying $\pi(\GD_{(s,t,\l)})=\ul D_{(t,\l)}$. 
 Notice that 
 for each $(t,\l)\in \Sigma_{\pi(z)}$ there is a unique $s=s(t,\l) \in S^1$ such 
 that $z\in \del \GD_{(s,t,\l)}$. 
 Actually if we write $z=[ ce^{h(\y)}: c e^{h(\y)}\y : \bar\y : 1]$ using 
 $\y\in\CP^1$ and $c\in \U(1)$, then such $s\in S^1$ is characterized by the 
 equation $e^{2H_++H_0+is}\Phi=ce^{h(\y)}$, 
 so $s$ is unique, and $s=s(t,\l)$ is smooth. 
 Hence we obtain a smooth section 
 $\Sigma_{\pi(z)} \to \cM : (t,\l)\to (s(t,\l),t,\l)$ 
 of which the image is $\GS_z\cap \cM$. 
 Hence $\GS_z\cap \cM$ is diffeomorphic to 
 $\Sigma_{\pi(z)}\cong \R\times S^1$. 
 
 Now we notice the disks $\{\GD_\x\}_{\x\in (S_+\sqcup S_-)}$. 
 Obviously, there are just two disks in this family satisfying 
 $z\in \del \GD_\x$. Hence $\GS_z$ is the natural compactification of 
 $\R\times S^1$ with extra two points, so $\GS_z$ is homeomorphic to $S^2$. 

 Finally we check the smoothness of $\GS_z$. 
 Let us put $\widetilde\GS_z:=\hat\Gf^{-1}(z)$.  
 By the similar argument as above, $\widetilde\GS_z\cap \cZ_\R$ is an 
 embedded $\R\times S^1$ in $\cZ_\R$, 
 and $\widetilde\GS_z$ is the natural 
 compactification of $\widetilde\GS_z\cap \cZ_\R$ with extra two points. 
 Recall that $\widetilde\GS_z\cap\cZ_\R=\Gf^{-1}(z)$ is an integral surface 
 of the distribution $\cD|_{\cZ_\R}$ by Proposition \ref{Prop:deformed_Gf}. 
 Then, the whole of $\widetilde\GS_z$ is an integral surface of $\cD$ 
 by the continuity. 
 Since the distribution $\cD$ is smooth, $\widetilde\GS_z$ is smoothly embedded 
 $S^2$. Hence $\GS_z=p(\widetilde\GS_z)$ is also a smoothly embedded $S^2$. 
\end{proof}


\begin{Prop} \label{Prop:extention}
 The map $\hat\Gf:(\hat\cZ_+,\hat\cZ_\R)\to (\CP^3,P_h)$ 
 satisfies the following conditions: 
 \begin{enumerate}
  \item each fiber of $\hat\Gf:\hat\cZ_\R \to P_h$ is an integral surface of 
    the distribution $\cD$, 
  \item $\hat\Gf: (\hat\cZ_+\backslash\hat\cZ_\R) \to 
     (\CP^3\backslash P_h)$ is biholomorphic,  
 \end{enumerate}
\end{Prop}
\begin{proof}
 We already showed the statement {\it 1} in the proof of Proposition 
 \ref{Prop:closedness_of_a_surf}. 
 The holomorphicity of $\hat\Gf : (\hat\cZ_+\backslash\hat\cZ_\R) \to 
   (\CP^3\backslash P_h)$  is deduced from {\it 1} 
 and the fiberwise holomorphicity of $\hat\Gf$. 
 Further, $\hat\Gf : (\hat\cZ_+\backslash\hat\cZ_\R) \to 
     (\CP^3\backslash P_h)$ is bijective by {\it 3} of 
 Proposition \ref{Prop:properties_of_disks}, hence {\it 2} follows. 
\end{proof}

\begin{Cor} \label{Cor:Zollfreiness}
 Let $(V,A)$ be an admissible monopole. Then the anti-self-dual metric on 
 $S^2\times S^2$ induced by $(V,A)$ is Zollfrei.  
\end{Cor}

\begin{proof} 
 As proved in \cite{bib:LM05} (Theorem 5.14), 
 an anti-self-dual 4-manifold $(S^2\times S^2,[g])$ is Zollfrei  
 if and only if every $\a$-surface is an embedded $S^2\subset S^2 \times S^2$. 
 (Here we are taking the opposite orientation to \cite{bib:LM05}.) 
 In our situation, every $\a$-surface is given as the image of 
 an integral surface of $\cD$ by $p:\hat\cZ_\R\to S^2\times S^2$. 
 Hence every $\a$-surface is written as 
 $\GS_z=\hat p(\hat \Gf^{-1}(z))$ for some $z\in P_h$ 
 by Proposition \ref{Prop:extention}. 
 Since $\GS_z$ is an embedded $S^2 \subset S^2\times S^2$ 
 by Proposition \ref{Prop:closedness_of_a_surf}, 
 the statement follows. 
\end{proof}

The proof of the main theorem (Theorem \ref{Thm:correspondence}) 
is already finished. 
Actually, the Zollfrei condition of the considering metric is 
proved in Corollary \ref{Cor:Zollfreiness}, and 
Proposition \ref{Prop:properties_of_disks} and \ref{Prop:extention} 
mean that the pair $(\CP^3,P_h)$ is the very LeBrun-Mason twistor space.

\section{Concluding remarks}

\paragraph{Regularity}
In this article, we assumed the smoothness of functions, embeddings and so on. 
In the previous articles \cite{bib:LM02,bib:LM05,bib:Nakata09}, however, 
we can construct the twistor correspondences of {\it low regularities}. 
Similar to these previous works, 
the argument in this article should be strengthened to that of low regularities. 
Actually the integral transforms $R$ and $Q$ are defined even for non-differentiable 
functions, and hyperbolic partial differential equations 
admit solutions of low regularities or distribution solutions in general. 
Thus the notion of a self-dual Zollfrei metric might be generalized to, 
for example, non-differentiable class. 
In fact, infinitely many examples of `self-dual Zollfrei metrics 
with singularity' are already obtained \cite{bib:Nakata07I}.

\paragraph{Degeneration}
We introduced the notion of admissible monopoles in 
Section \ref{Sect:Monopole_equation}, 
and showed that the corresponding admissible deformations 
$\RP^3$ in $\CP^3$ has nice properties and the 
LeBrun-Mason correspondence works well (Theorem \ref{Thm:correspondence}). 
On the other hand, in the non-admissible case, 
the deformation of $\RP^3$ in $\CP^3$ 
has an unexpected property (Proposition \ref{Prop:non-admissible_case}). 
Even in the non-admissible case, however, 
we can get the family of holomorphic disks 
parametrized by $S^2\times S^2$ 
(Proposition \ref{Prop:deformed_holomorphic_discs}). 
Then the natural question is: 
\begin{itemize}
\item Is there any natural structure on the 
parameter space of the holomorphic disks for the non-admissible case?  
\end{itemize}
In particular, it would be interesting to study the process of the degeneration 
which occurs in the deformation from an admissible case to a non-admissible case.

\paragraph{Deformation of $S^3_1$}
The argument in this article is based on the 
identification between the de Sitter space $S^3_1$ and the space of 
oriented small circles on the two sphere $\bS^2$, which is 
arisen from the LeBrun-Mason correspondence for Einstein-Weyl structures
\cite{bib:LM08,bib:Nakata09}. 
By the result in \cite{bib:LM08,bib:Nakata09}, 
if we deform the twistor space 
from $(W,W_\R)=(\CP^1\times\CP^1,\CP^1)$ to $(W,P)$, 
we obtain an 
Einstein-Weyl structure on $\R\times S^2$ of indefinite signature. 
In this construction, 
$\R\times S^2$ is identified with the space of oriented circles embedded in 
$P\simeq\CP^1$. 
So it is natural to expect the generalization of our story to such 
deformed situations. 
If it is successful, we will obtain various significant objects: 
general solutions of the wave equations on $\R\times S^2$, 
descriptions of more general self-dual Zollfrei metrics, 
its LeBrun-Mason twistor spaces, and so on.

\vspace{10mm}

\appendix

\noindent
{\Large \bf Appendix}

\vspace{-3mm}
\section{The bijectivity of $\cQ$}

We give a proof of the bijectivity of the transform 
$\cQ: C^\infty_{\text{odd}}(\bS^2) \to C^\infty_{\text{odd}}(\bS^2)$ 
by a similar method as Guillemin's \cite{bib:Guillemin76}. 
Let $\hat\scrH^k$ be the space of homogeneous harmonic 
polynomials of degree $k$ on $\R^3$ and let 
$\scrH^k=\{P|_{\bS^2} \in L^2(\bS^2) \mid P\in\hat\scrH\}$.  
We notice the following fact. 

\begin{Thm} \label{Thm:decomp.of_L2S2}
 The group $\SO(3)$ acts irreducibly on $\scrH^k$ and 
 the representations on $\scrH^k$ and $\scrH^l$ are inequivalent if $k\neq l$ . 
 Moreover, we have the decomposition 
  \begin{equation} \label{eq:decomp_of_L^2(S^2)}
   L^2(\bS^2)\cong \oplus_k \scrH^k
  \end{equation}
 as a direct sum of Hilbert spaces. 
\end{Thm}

Since $\cQ$ maps $L^2(\bS^2)$ to itself and commutes with the 
$SO(3)$-action, so $\cQ$ is diagonalized with respect to the decomposition 
\eqref{eq:decomp_of_L^2(S^2)}. Let us denote the eigenvalues of 
$\cQ$ on $\scrH^k$ by $c(k)\in \R$, that is, 
$$ \cQ h= c(k) \cdot h \qquad \text{for} \quad h\in\scrH^k. $$

\begin{Prop} \label{Prop:A2}
 $$ c(k) = \begin{cases}
   1  & k=0, \\ 
   0  & k=2m \quad (m=1,2,\cdots), \\
   (-1)^m \, \frac{4\pi}{2m+1} \cdot 
  \frac{1\cdot 3\cdot 5\cdots (2m+1)}{2\cdot 4\cdot 6\cdots (2m+2)} 
    \ \ & k=2m+1 \ \ (m=0,1,\cdots). 
  \end{cases} $$
\end{Prop}
\begin{proof}
 Since $\cQ(1)=1$ by definition, we obtain $c(0)=1$. 
 On the other hand, 
 since $C^\infty_{\text{even}\, *}(\bS^2)$ is annihilated by $\cQ$, 
 we obtain $c(2m)=0$ for $m>0$. 
 
 Suppose $k=2m+1$. 
 Let us choose a harmonic polynomial $P(x,y,z)\in \hat\scrH^k$ so that 
 it does not depend on the $z$ variable. Then $P$ is written as 
 $$ P(x,y,z)=a_{2m+1}x^{2m+1} + a_{2m}x^{2m}y+\cdots +a_0 y^{2m+1}.$$
 Since $P$ is harmonic, the equation $(\del_x^2+\del_y^2)P=0$ holds. 
 Hence we obtain 
 \begin{equation} \label{eq:A2}
  a_{2j+1}= - \frac{(2m-2j+2)(2m-2j+1)}{(2j+1)\cdot 2j}\, a_{2j-1} 
  \qquad (j=1,2,\cdots,m)   
 \end{equation}
 \begin{equation} \label{eq:A3}
  \text{or} \qquad a_{2m+1}= (-1)^m \, 
  \frac{2m\cdot(2m-1)\cdots 2\cdot 1}{(2m+1)\cdot 2m\cdots 3\cdot 2}\, a_1
  =  \frac{(-1)^m}{2m+1}\, a_1.  \hspace{22mm}
 \end{equation}
 Now we have $(\cQ P)(1,0,0)=c(2m+1)P(1,0,0)= c(2m+1) a_{2m+1}$. 
 On the other hand, by definition, 
 $$ (\cQ P)(1,0,0)=\int_\Omega P(x,y)\w_{\bS^2} \qquad \text{where} \quad
  \Omega=\{(x,y,z)\in \bS^2 \mid x>0 \}. $$
 Let us use the coordinate $(\c,\vp)$ so that 
 $$ (x,y,z)=(\sin\c\cos\vp,\, \sin\c\sin\vp,\, \cos\c), \qquad 
  \Omega=\left\{ 0\le\c\le\pi, 
            -\frac{\pi}{2}\le\vp\le\frac{\pi}{2}\right\}. $$	
 Then, by $\w_{\bS^2}= \sin\c\, d\c d\vp$,   
 $$ \begin{aligned}  (\cQ P)(1,0,0) 
   &=\int_\Omega P(\sin\c\cos\vp,\sin\c\sin\vp) \sin\c \, d\c d\vp \\[1ex]
   &= \sum_{l=0}^{2m+1} a_l \( \int_0^\pi (\sin\c)^{2m+2} d\c \)
       \(\int_{-\frac{\pi}{2}}^{\frac{\pi}{2}} 
	   (\cos\vp)^l (\sin\vp)^{2m-l+1} d\vp\) \\ 
   &= 2 \( \int_0^\pi (\sin\c)^{2m+2} d\c \) \sum_{j=0}^m a_{2j+1}  
       \(\int_0^{\frac{\pi}{2}} (\cos\vp)^{2j+1} (\sin\vp)^{2m-2j} d\vp\). 
	 \end{aligned} $$
 By a usual trick, which is also explained in \cite{bib:Guillemin76}, we obtain 
 \begin{equation} \label{eq:A4}
   \int_0^\pi (\sin\c)^{2m+2} d\c = 2\pi \, 
	\frac{1\cdot 3 \cdot 5 \cdots (2m+1)}{2\cdot 4\cdot 6\cdots (2m+2)}. 
 \end{equation}
 On the other hand, if we put 
 \begin{equation}
  B(j)=\int_0^{\frac{\pi}{2}}  (\cos\vp)^{2j+1} (\sin\vp)^{2m-2j} d\vp, 
 \end{equation}
 then for $j>0$ 
 $$ \begin{aligned} B(j) &=\left[  (\cos\vp)^{2j}\cdot 
              \frac{(\sin\vp)^{2m-2j+1}}{2m-2j+1} \right]_0^{\frac{\pi}{2}}
           - \int_0^{\frac{\pi}{2}} \pd{\vp}((\cos\vp)^{2j}) \cdot 
	          \frac{(\sin\vp)^{2m-2j+1}}{2m-2j+1} d\vp \\
		&= \frac{2j}{2m-2j+1} B(j-1). \end{aligned} $$
 Hence, combining with \eqref{eq:A2}, we obtain 
 $$ a_{2j+1}B(j) = - \frac{2m-2j+2}{2j+1}\, a_{2j-1}B(j-1) \hspace{21.6mm}$$ 
 \begin{equation} \label{eq:A6}
  \text{or} \qquad a_{2j+1}B(j) 
   = \frac{2m-2j}{2m+1} \, a_{2j+1}B(j) 
    - \frac{2m-2j-2}{2m+1}\, a_{2j-1}B(j-1). 
 \end{equation}
 If we take a sum of \eqref{eq:A6} for $j=1,2,\cdots,m$, then we obtain 
 $$ \sum_{j=1}^m a_{2j+1}B(j)= - \frac{2m}{2m+1}\, a_1 B(0). $$ 
 Thus 
 \begin{equation} \label{eq:A7}
   \sum_{j=0}^m a_{2j+1}B(j)= \frac{1}{2m+1}\, a_1 B(0)
     = \frac{(-1)^m}{2m+1} a_{2m+1}. 
 \end{equation}
 By \eqref{eq:A4} and \eqref{eq:A7},  
 $$ (\cQ P)(1,0,0)= c(2m+1)a_{2m+1} = 
  (-1)^m \, \frac{4\pi}{2m+1} \cdot 
  \frac{1\cdot 3\cdot 5\cdots (2m+1)}{2\cdot 4\cdot 6\cdots (2m+2)} 
  a_{2m+1}. $$
 Since $a_{2m+1}\neq 0$, we obtain the required formula. 
\end{proof}

We denote the degree $s$ Sobolev space over $\bS^2$ by $H^s$, 
and put $H^s_{\text{odd}}:=H^s\cap L^2_{\text{odd}}(\bS^2)$. 
Let us define a norm on $H^s$ by 
$$ |h|_s = k^s |h|_{L^2} \qquad \text{for} \quad h\in\scrH^k,$$
then, as explained in \cite{bib:Guillemin76}, this norm is equivalent 
to the usual $H^s$ norm. 
\begin{Prop} \label{Prop:A3}
 There exists a constant $c>1$ independent of $s$ such that 
 $$ \frac{1}{c} |h|_s \le |\cQ h|_{s+\frac{3}{2}} \le c |h|_s $$
 for all $s$ and $h\in H^s_{\text{odd}}$. 
 Hence $\cQ$ defines a bijection 
 $H^s_{\text{odd}}\to H^{s+\frac{3}{2}}_{\text{odd}}$. 
\end{Prop}
\begin{proof}
 Similarly to \cite{bib:Guillemin76}, we notice the formula 
 $$ \pi^{\frac{1}{2}}= \lim_{k\to\infty} k^{-\frac{1}{2}}\, 
   \frac{2\cdot 4 \cdot 6\cdots 2m}{1\cdot 3\cdot 5\cdots (2m-1)}. $$
 By Proposition \ref{Prop:A2}, we get 
 $$ c(2m-1) \sim (-1)^k 2 \pi^{\frac{1}{2}} k^{-\frac{3}{2}}. $$
 So the statement follows. 
\end{proof}

By Proposition \ref{Prop:A3} and the Sobolev's embedding theorem, 
we obtain the following. 
\begin{Thm}
 The transform $\cQ: C^\infty_\text{odd}(\bS^2) 
  \to C^\infty_\text{odd}(\bS^2)$ is bijective. 
\end{Thm}

\section{Formulas on $S^3_1$} \label{App:Formulas}
In Section \ref{Sect:Standard_model} equation \eqref{eq:onf_S^3_1}, 
we introduced a local orthonormal frame $\{\ul E_1,\ul E_2,\ul E_3\}$ 
of the tangent bundle $TS^3_1$ on the open set 
$\ul U=\{(t,\l)\mid \l\neq\infty \} \subset S^3_1$. 
Here we show several formulas concerning this frame. 
All these formulas are deduced by direct calculations.

The connection form $\w$ of the Levi-Civita connection for $g_{S^3_1}$ and its curvature form $K$ 
is 
\begin{equation}
 \w= \begin{pmatrix} 0&\w^1_2&\w^1_3 \\ \w^1_2&0&\w^2_3 \\
   \w^1_3&-\w^2_3&0 \end{pmatrix}
 \qquad
 \begin{aligned}
  \ul \w^1_2 &= \tanh t \, \ul E^2, \qquad
  \ul \w^1_3 = \tanh t \, \ul E^3, \\
  \ul \w^2_3 &= \frac{1}{\cosh t}\( -\Imag\l\, \ul E^2 + \Real\l\, \ul E^3 \), 
 \end{aligned}
\end{equation}
\begin{equation}
 K= \begin{pmatrix} 0 & \ul E^1\wedge \ul E^2 & \ul E^1\wedge \ul E^3 \\ 
      \ul E^1\wedge \ul E^2 & 0 & \ul E^2\wedge \ul E^3 \\
      \ul E^1\wedge \ul E^3& \ul E^3\wedge \ul E^2&0  \end{pmatrix}.
\end{equation}
Let $\ul\Gm_1(\z)$ and $\ul\Gm_2(\z)$ be the vector fields 
on $\ul U$ defined by \eqref{eq:EW_Gm1,2}.  
By the lifting formula \eqref{eq:EW_lift}, the tautological lifts 
$\tilde{\ul\Gm}_1$ and $\tilde{\ul\Gm}_2$ on $\cW_\R$ (or on $\cW_+$)
are written as 
\begin{equation}
 \left\{ \begin{aligned}
  \tilde{\ul\Gm}_1&= \ul\Gm_1 + \g_1 \del_\z, \\
  \tilde{\ul\Gm}_2&= \ul\Gm_2 + \g_2 \del_\z, 
 \end{aligned} \right. 
 \qquad
 \begin{aligned}
  \g_1&= \Psi\cdot (-\Imag\l + \z \Real\l - \z \sinh t), \\ 
  \g_2&= \Psi\cdot (-\z \Imag\l -\Real\l - \sinh t), \\ 
 \end{aligned}
 \qquad \Psi := \frac{1+\z^2}{2\cosh t}. 
\end{equation}
If we change the fiber coordinate by $\z=i\, \frac{1-\w}{1+\w}$, 
we obtain 
\begin{equation} \label{eq:Gm_in_omega}
 \left\{ \begin{aligned}
 (1+\w)\tilde{\ul \Gm}_1 &= l_1 +\w\, \bar l_1 
    + (\d_1-\w\, \bar \d_1)\, \w \del_\w, \\
 (1+\w)\tilde{\ul \Gm}_2 &= l_2 +\w\, \bar l_2 
    + (\d_2-\w\, \bar \d_2)\, \w \del_\w, 
 \end{aligned} \right. 
\end{equation}
$$ \text{where} \qquad \left\{
   \begin{aligned}
    l_1 &= -\ul E_1+ \ul E_2 +i \ul E_3, \quad 
    \d_1 = -\frac{\l}{\cosh t}+\tanh t, \\ 
	l_2 &=i\ul E_1+ i \ul E_2 - \ul E_3, \quad \ 
	\d_2 = -\frac{i\l}{\cosh t}-i\tanh t.  
 \end{aligned} \right. $$

\vspace{5mm}
\paragraph{Acknowledgement.}
The author would like to thank Hiroyuki Kamada for helpful conversations 
and for providing many useful references. He would also like to thank Mikio Furuta 
for helpful discussions.



\vspace{13mm}
\noindent
\small
\begin{tabular}{l}
Department of Mathematics \\
Graduate School of Science and Engineering \\
Tokyo Institute of Technology \\
2-12-1, O-okayama, Meguro, 152-8551, JAPAN \\
{\tt {nakata@math.titech.ac.jp}}
\end{tabular}

\end{document}